\documentclass[12pt]{article}
\usepackage{amsmath, amssymb, amsthm}
\usepackage[all]{xypic}
\usepackage{fullpage}
\newtheorem{thm}[equation]{Theorem}
\newtheorem{prop}[equation]{Proposition} 
\newtheorem{cor}[equation]{Corollary}
\newtheorem{lemma}[equation]{Lemma}

\theoremstyle{definition}
\newtheorem{defn}[equation]{Definition}

\theoremstyle{remark}

\newtheorem{ntn}[equation]{Notation}

\newtheorem{rem}[equation]{Remark}

\makeatletter
\renewcommand{\subsection}{\@startsection{subsection}{2}{0pt}{-3ex
plus -1ex minus -0.2ex}{-2mm plus -0pt minus
-2pt}{\normalfont\bfseries}} \makeatother

\numberwithin{equation}{subsection}

\newdir^{ (}{{}*!/-5pt/@^{(}}

\newcommand{\iso}{{\;\stackrel{_\sim}{\to}\;}}

\newcommand{\beq}{\begin{equation}\label}
\newcommand{\eeq}{\end{equation}}
\def\ccirc{{{}_{\,{}^{^\circ}}}}

\newcommand{\onto}{\twoheadrightarrow}
\newcommand{\tooo}{{\;{-\!\!\!-\!\!\!-\!\!\!-\!\!\!\longrightarrow}\;}}

\newcommand{\pr}{\operatorname{pr}}

\newcommand{\ort}{\text{or}}
\newcommand{\mfg}{\mathfrak g}

\newcommand{\Hom}{\text{Hom}}

\def\R{\mathbb{R}}

\def\k{\mathbf{k}}

\def\lb{{\mathbf{l}}}
\def\ppi{\mathrm{pr}}
\def\o{\otimes}
\def\h{\mathrm{h}}

\def\dq{\overline{Q}}

\def\Id{\mathrm{Id}}

\def\moy{{\mathrm{Moyal}}}
\def\Rep{{\mathrm{Rep}}}
\def\d{{\mathbf{l}}}
\def\Sym{{\text{Sym\ }}}
\def\Z{{\mathbb Z}}
\def\Q{{\mathbb Q}}
\def\tr{{\text{tr}}}
\def\lnpp{L_{n,+}}
\def\lnp1p{L_{n+1,+}}
\def\lnqpp{L_{nQ,+}}
\def\lnqp1p{L_{(n+1)Q,+}}
\def\gl{\mathfrak{gl}}
\def\sp{\mathfrak{sp}}

\begin{document}
\title{Moyal quantization and stable homology of necklace Lie algebras}
\author{Victor Ginzburg and Travis Schedler}
\maketitle




\begin{abstract} 
We compute the stable homology of  necklace Lie
algebras  associated with quivers and give a construction of stable
 homology classes from certain $A_\infty$-categories.  Our construction
 is a generalization of the construction of homology 
classes of moduli spaces of
 curves due to M. Kontsevich.

In the second part of the paper we produce a  Moyal-type 
 quantization of the symmetric algebra of a  necklace Lie
algebra.
The resulting quantized algebra has  natural  representations
in the usual Moyal quantization of polynomial algebras.
\end{abstract}


This paper consists of two parts, essentially independent of each other.
The results of the first part, to be outlined in Sect.
\ref{intro1}, are concerned with stable homology of
necklace Lie algebras. In the second part, we will be concerned
with a natural quantization of necklace Lie algebras.
The results of this part are outlined in Sect.
\ref{intro2} below.
 
\section{Stable homology of necklace Lie algebras.}\label{intro1}
\subsection{The graph complex and Lie algebra homology.}
In \cite{K}, Kontsevich gives a construction for the stable homology of
Lie algebras associated to commutative, Lie, and associative operads
(generalized to arbitrary operads in \cite{CV}).  A key idea of
Kontsevich was to interpret the chain complex
involved in the computation  of 
Lie algebra  homology in question as a certain
\textsl{graph complex}. Furthermore,
in the associative and Lie cases,  Kontsevich related 
the homology of the graph complex with
the cohomology of the coarse moduli space of smooth algebraic curves of
genus $g$ with $n$ punctures, and the space of outer automorphisms of a
free group with $n$ punctures, respectively.

The Lie algebra corresponding to the associative operad is defined in
\cite{K} as follows.  Throughout, fix a field $k$ of characteristic
zero.  For any $n$, let $P_n$ be the free associative (noncommutative)
$k$-algebra with generators $x_1, x_2, \ldots, x_n, y_1, \ldots, y_n$.
Let $L_n$ be the sub-Lie algebra of derivations of $P_n$ which kill
the element $\sum_{i=1}^n [x_i, y_i] \in P_n$, which can be interpreted
as ``Hamiltonian vector fields'', and let $\lnpp \subset
L_n$ be the subspace spanned by derivations of nonnegative degree
(where a derivation has degree $d$ if it sends homogeneous polynomials
of degree $e$ to degree $d+e$).  The latter can be interpreted as
``Hamiltonian vector fields that fix the origin'': these do not include the
$\frac{d}{dx_i} \in L_n$.  There are natural inclusions $\lnpp
\subset \lnp1p, P_n \subset P_{n+1}$, which induces maps on Lie
algebra homology. It makes sense therefore to consider the
\textsl{stable homology}, $\underset{\rightarrow}{\lim} H_*(\lnpp) =
H_*(L_\infty)$, where $L_\infty = \bigcup_{n=1}^\infty \lnpp$.  (We
omit the ``$+$'' since we will not consider the stable version without
a plus).  A standard argument based on an inclusion $L_\infty \oplus
L_\infty \rightarrow L_\infty$ shows that $H_*(L_\infty)$ has a
natural structure of a graded cocommutative Hopf algebra.  Write
$PH_*(L_\infty)=\bigoplus_{k\geq 0}PH_k(L_\infty)$ for the
corresponding graded vector space of primitive elements, and let
$H^*_c$ stand for cohomology with compact support.

In
\cite{K}, Kontsevich
established an isomorphism
\begin{equation}\label{kon1}
PH_k(L_\infty) \cong PH_k(\sp(2\infty)) \oplus \bigoplus_{m > 0, 2-2g-m < 0}
 H_c^{2g-2+m+k} (M^{\text{comb}}_{g,m}/\Sigma_m, k),
\end{equation}
where $M^{\text{comb}}_{g,m}$ is an orbicell
complex whose associated chain complex
is known as the \textbf{ribbon graph complex}.
This is a chain complex, $RG^{g,m}_\cdot$,
whose terms are vector spaces with bases labelled by connected ribbon graphs
(whose vertices have valence $\geq 3$) with two fixed combinatorial
invariants: $g=${\em genus}, and $m=${\em number of punctures}. Passing
to {\em primitive} homology above corresponded to restricting here to
{\em connected} ribbon graphs. Taking
quotient by the action of the symmetric group $\Sigma_m$ in the right
hand side of \eqref{kon1} amounts to forgetting the order of the
punctures.  The differential is defined via edge contractions.

The notation of $M^{\text{comb}}_{g,m}$
 is due to a homeomorphism
 $M^\text{comb}_{g,m} \cong \mathcal
 M_{g,m} \times \R^{m}$, where  $\mathcal M_{g,m}$ is
the coarse moduli space  of
smooth complex algebraic curves of genus $g$ with $m$ punctures.
With this homeomorphism in mind, passing to 
primitive homology corresponds to considering
connected curves only.

\begin{rem} In \cite{K} the result is stated as an isomorphism with
$H^{4g-4+2m-k}$ of $\mathcal M_{g,m}$, presumably using Poincar\'e
 duality for the orbifold $M^\text{comb}_{g,m}$, cf. Remark \ref{qpd},
 and the homeomorphism $M^\text{comb}_{g,m} \cong \mathcal M_{g,m}
 \times \R^m$.
\end{rem}
\subsection{A quiver analogue.}
The story in the previous subsection has a generalization to quivers.
Namely, the free associative algebra $P_n$ can be viewed as the path
algebra of the quiver with one vertex and $2n$ loops.  Further, the
algebra of derivations killing $\sum_{i} [x_i, y_i]$ is the
so-called {\em necklace Lie algebra} for this quiver.  One can more
generally associate a necklace Lie algebra to any quiver $Q$ as
follows (from \cite{G} and \cite{BLB}).  First take the double quiver
obtained by adding a reverse edge $e^*\in\dq$ for each edge $e\in
Q$. Let $P_Q$ be the path algebra of $\dq$.  Then we can consider the
Lie algebra $L_Q$ of derivations of $P_Q$ killing the element
$\sum_{e \in Q} [e, e^*] \in P_Q$. Let $L_{Q,+} \subset L_Q$ restrict
to the span of derivations of nonnegative degree (as in the case $L_n$: they
do not decrease the degree of homogeneous elements). As before, $L_{Q}$
and $L_{Q,+}$ can be interpreted as the Lie algebra of Hamiltonian
vector fields and the subalgebra which fixes the origin. (See also
Definition \ref{lqigd} about the grading on $L_Q, L_{Q,+}$).

To form a stable version, we consider, for any quiver $Q$, the quiver
$nQ$ with the same vertex set $I$ as $Q$, obtained by taking $n$ copies
of each edge of $Q$.  Then, for any $n$, we can
consider $P_{nQ}, \lnqpp$.  Once again, we have natural inclusions
$\lnqpp \subset \lnqp1p$ and $P_{nQ} \subset P_{(n+1)Q}$.  So,
analogously to the case where $Q$ is a quiver with just one vertex and
one edge, one can consider the stable homology
$\underset{\rightarrow}{\lim} H_*(\lnqpp) = H_*(L_{\infty Q})$, where
$L_{\infty Q} := \bigcup_{n=1}^\infty \lnqpp$.  As before, this is a
Hopf algebra, and one can consider the primitive part $PH_*(L_{\infty
Q})$ (in other words, one can restrict to the connected graphs that will
appear, and avoid the consideration of disconnected ones).

In fact, the stable Lie algebra $L_{\infty Q}$ does not depend on the
multiplicities of edges in $Q$, except whether the multiplicity of edges
from a vertex $i$ to $j$ is zero or nonzero.  This is because the stable
quiver $\infty Q$ is simply a quiver with edge multiplicities equal to
$0$ or $\infty$.  As a result, the stable homology depends only on the
adjacency as well.

So from now on, we replace the quiver $Q$ by an
undirected graph $G$ with vertex set $I$, such that the
multiplicity of edges from $i$ to $j$ is either one or zero.
(Of
course, $G$ may still have loops, but only one or zero at each vertex.)
Write $v
\sim w$ if $v$ and $w$ are adjacent in $G$.    We call the \textbf{stable cohomology of the necklace Lie algebra for $G$} the stable cohomology as constructed above for any quiver whose double has the same adjacencies as $G$ (without multiplicity).

We then describe this cohomology similarly to Kontsevich's description in the
case of one vertex.
We will define a topological space $M^{\text{comb}}_{g,m,G,X}$,
which is related to
 $M^{\text{comb}}_{g,m}$. Let $I^{(m)}$ denote the set of unordered $m$-tuples
of
elements of $I$. Essentially, the space  $M^{\text{comb}}_{g,m,G,X}$ is
obtained by taking
 isomorphism classes of ribbon graphs with a metric and $n$ numbered
faces (punctures), all of whose vertices have
 valence $\geq 3$, together with an additional labeling of the faces 
by the terms in $X \subset I^{(m)}$ such
 that the only adjacent faces (meaning that there is an edge which meets
 the two faces) are those whose $I$-labels are adjacent in $G$.  This includes
 self-adjacency: 
an edge can meet two boundary
 components with the same label (or a single boundary component) only if
the corresponding vertex of $G$ has a loop.

We then mod by the symmetric group $\Sigma_m$ as before, to get an
orbicell complex whose associated compactly-supported cohomology is
the stable cohomology of the necklace Lie algebra for the graph $G$.


Before proceeding, we give some comments about this definition:
\begin{rem} The quotient $M^{\text{comb}}_{g,m,G,X}/\Sigma_m$ is a
subquotient of $M^{\text{comb}}_{g,m}$ (the latter is the usual,
 non-quiver version).  Namely, let's replace each $X$ with a fixed lift
 $\tilde X \in I^n$: for example, order $I$ and let the fixed lift be
 those that are in order lexicographically ($\tilde X =
 (a_1,\ldots,a_m)$ where $a_i \leq a_j$ in the ordering $\leq$ on $I$).
 Then we consider the subcomplex of $M^{\text{comb}}_{g,m}$ of ribbon
 graphs such that if we label the faces by $\tilde X$ (which says for
 each puncture what vertex it should be labeled by), then each pair of
 adjacent faces have adjacent vertex labels.  Then, we can quotient by
 the subgroup of $\Sigma_{m}$ which stabilizes $\tilde X$.
\end{rem}

\begin{rem} 
We can describe the construction of $M^{\text{comb}}_{g,m,G,X}$ using
the dual of the ribbon graph as follows: Giving a labeling of the
faces of a ribbon graph $\Gamma$ by $I$ such that adjacent faces have
adjacent labels is the same as giving a morphism from the dual of
$\Gamma$ to the graph $G$: more precisely, forget the ribbon graph
structure to take the underlying undirected graph of the dual of
$\Gamma$, then give a morphism of graphs of this to $G$.

The trivalence condition then becomes the condition that each face
have at least three edges.  Then we get the complex of isomorphism
classes of ribbon graphs $\Gamma$ with labeled vertices and a metric
(on edges), all of whose faces have at least three edges, together
with a map $\eta$ preserving adjacency to $G$. [When we say
isomorphism classes, we mean where isomorphisms are isomorphisms of
ribbon graphs with metric, preserving the labeling $1,\ldots,m$,
commuting with the map $\eta$ to $G$.]

The differential on the dual complex is really simple: it just involving
 summing over which edge to delete (not contract, simply delete), with
 sign.  We simply don't allow removing an edge that would disconnect the
 graph.  

To get the complex $M^{\text{comb}}_{g,m,G,X}/\Sigma_m$, we can forget
the labeling requirement of the vertices of the ribbon graph: then we
take isomorphism classes of unlabeled ribbon graphs $\Gamma$ together with
a map $\eta: \Gamma \rightarrow G$.

Of course, the dual point of view is formally equivalent to the original
one, so we will not use this interpretation.
\end{rem}

The problem with defining 
$\mathcal M_{g,m,G,X}$ (non-combinatorial version) is that for a point in the moduli space $\mathcal M_{g,m}$ without metric, it's not obvious whether it corresponds to a ribbon graph that can be labeled by $X$ or not.

On the other hand, note that there is a more conceptual definition of
the space $\mathcal M_{g,m}$
(kindly explained to us by K. Costello) as the geometric realization
of a category whose objects are ribbon graphs and whose morphisms are
quotients by sub-forests (unions of trees with no common vertices).
Similarly, one could define the space $\mathcal M_{g,m,G,X}$ as the
geometric realization of the category whose objects are $X$-labeled
ribbon graphs.  We can think of this $\mathcal M_{g,m,G,X}$ as the
``moduli space of $(G, X)$-labeled ribbon graphs of genus $g$'' or
``moduli space of $(G, X)$-labeled surfaces of genus $g$''.  Then, one
should be able to get a similar formula to the following Theorem which
uses $\mathcal M_{g,m,G,X}$ instead of $M^{\text{comb}}_{g,m,G,X}$ and
avoids metrics.

We now state our first main result.
\begin{thm} \label{sht} 
The primitive stable homology $PH_k(L_{\infty G})$ is given
 by
\begin{multline}
PH_k(L_{\infty G}) \cong \bigoplus_{v \sim v} PH_k(\sp(2 \infty)) \oplus
 \bigoplus_{v \sim w, v \neq w}  PH_k(\gl(\infty)) \\ \oplus
 \bigoplus_{X \in I^{(m)}, 2-2g-m < 0} H_c^{2g-2+m+k} 
(M^{\text{comb}}_{g,m,G,X}).
\end{multline}
\end{thm}

The above isomorphism actually comes from an isomorphism of some natural
cochain complexes which compute cohomology.

\begin{rem} \label{phkcr}
The $PH_k(\sp(2 \infty))$ and $PH_k(\gl(\infty))$ terms come from the
part which has vertices with valence $\leq 2$ (cf. Lemma \ref{ssl}),
and are easy to compute, for instance by viewing them as the polygons
in the graph complex.  The result is $PH_k(\sp(2 \infty)) = \Q$ if $k
\equiv 3 \pmod 4$ and $0$ otherwise, and the same for
$\gl(\infty)$.
\end{rem}

\begin{rem}\label{qpd} It would be interesting to find out in which cases one can
apply a Poincar\'e duality to this orbicell complex, which can always
be done in the case that the quiver has one vertex.  Then, one could
replace cohomology with compact support with regular cohomology.  For
example, with one vertex and at least one loop, the top cohomology
classes correspond to ribbon graphs where all vertices are trivalent
(note that in this case, one also knows that Poincar\'e duality
applies since $\mathcal M_{g,m}$ itself is an orbifold). However, in
the quiver case, there is not always a nonzero ribbon graph with all
vertices trivalent: for example, this does not exist if any of the
elements of $I$ appearing in $X$ are only vertices of even-sided
cyclic paths in $G$ (or if all of the odd-length cyclic paths involve
vertices not appearing in $X$).  However, it may still be possible to
apply Poincar\'e duality in such cases using a lower-dimensional top
class.
\end{rem}

\subsection{Stable homology classes from  $A_\infty$-algebras.} \label{ainfcs}
Kontsevich  showed, see \cite{KFD}, that any
cyclic finite-dimensional $A_\infty$ algebra with an inner product
 gives rise to a cycle in the 
ribbon graph complex. 
Specifically, let $A$ be a
finite-dimensional  and $\Z/2$-graded $A_\infty$ algebra with an even 
cyclic nondegenerate
symmetric bilinear form $\langle-,- \rangle: A\times A\to k$.
Thus, for each $n=1,2,\ldots,$
there is an $n$-ary operation $m_n: A^{\otimes n} \rightarrow A[n]$,
a graded map of parity  $|m_n| = 2-n$.  For simplicity, we assume that
$m_1 = 0$.
 Associated with $m_n$, is the pairing $\tilde m_{n}: A^{\otimes (n+1)}
\rightarrow k$ by $(x_1 \otimes x_2 \otimes \cdots x_{n+1}) = \langle
m_n(x_1 \otimes \cdots \otimes x_n), x_{n+1} \rangle)$.
The axioms of  a cyclic  $A_\infty$ algebra read
\begin{gather}
\sum_{k,\ell} (-1)^{\ell(d_1+\ldots+d_k)+(k+1)(\ell+1)} m_{n-\ell+1} (1^{\otimes k}
 \otimes m_\ell \otimes 1^{\otimes (n-\ell-k)}) = 0, \forall n
 \label{ainfid}\\ 
\tilde m_n(v_2 \o \cdots \o v_{n+1} \o v_1)  = (-1)^{n + (d_1) (d_2
 + d_3 + \ldots + d_{n+1})}\tilde m_n(v_1 \o v_2 \o
 \cdots \o v_{n+1}), \forall n\geq 1,
\end{gather}
where $d_i(v_1 \otimes \cdots \otimes v_n) = |v_i| \in \Z/2$, so that
the term $(-1)^{\ell(d_1 + \ldots + d_k) + (k+1)(\ell+1)}$ in
\eqref{ainfid} becomes $(-1)^{\ell(|v_1| + \ldots + |v_k|) +
(k+1)(\ell+1)}$ when the summand is applied to an element $v_1 \otimes
\cdots \otimes v_n$.

The pairing $\tilde m_n, n\geq 2$ may be viewed as an element of
 $(A^{\otimes(n+1)})^*,$ which is graded cyclically-symmetric and even
 for $n$ even and graded cyclically-antisymmetric and odd for $n$ odd.
 For $n=1,$ the inverse to the nondegenerate bilinear form $\langle-,-
 \rangle$ gives a graded symmetric even element $C \in
 A^{\otimes(2)}$.

To such a cyclic $A_\infty$ algebra $A$, Kontsevich associates a chain
in the ribbon graph complex, that is, a formal linear combination of
various oriented graphs with certain coefficients, called {\em weights}
(for the definition of orientation on graphs, see Section \ref{rgcss}).  To
define the weight corresponding to an oriented ribbon graph, observe that the
ribbon graph defines precisely a way to contract copies of the tensors
$\tilde m_n, C$: namely, we take the product over all vertices of
$\tilde m_n$ in some order, and over all edges of $C$, yielding an
element of $(A^{\otimes (2\#(E))}) \o (A^{\otimes (2\#(E))})$.  Then
we perform graded contractions to get an element of $k$.  Here, a
graded contraction means $\langle v_1 \o \cdots \o v_m, w_1 \o \cdots
\o w_m \rangle_{\theta} = \pm \langle v_1, w_{\theta(1)} \rangle
\cdots \langle v_m, w_{\theta(m)} \rangle$, where the sign $\pm$ is
obtained by applying the braiding $w_i \o w_j \mapsto (-1)^{|w_i|
|w_j|} w_j \o w_i$ to the right component many times to obtain $\pm
w_{\theta(m)} \o w_{\theta(m-1)} \o \cdots \o w_{\theta(1)}$.  (This
could be more naturally defined by using the language of braided
tensor categories.)  The order that the vertices and edges were placed
in only creates an ambiguity of sign, which is removed by using the
orientation on the graph.

It turns out that equation \eqref{ainfid} guarantees that
the chain constructed via the above  procedure is a
{\em cycle} in the ribbon graph complex.  We refer to
 \cite{PS},\cite{Pk} for more details.

\subsection{Quiver generalization and  $A_\infty$-categories.} \label{ainfcs2}
In this paper, we extend the above construction by Kontsevich to the
quiver setting.  In more detail, let $G$ be a graph with edge
multiplicities $0$ or $1$. As explained earlier, such a graph gives
rise to a necklace Lie algebra $L_{\infty G}$. We show, generalizing
Kontsevich's construction, that any cyclic $A_\infty$-{\em category}
with inner product gives rise to a cycle in the chain complex for
$M^{\text{comb}}_{g,m,G,X}$. The objects of the $A_\infty$-category in
question correspond to the vertices of the graph $G$ and morphisms are
generated by the morphisms from edges of $G$.

In more detail, we consider a structure which assigns to any edge of
$G$ with endpoints $i, j \in I$ with $i \neq j$ two
finite-dimensional, $\Z/2$-graded vector spaces $V_{ij}, V_{ji}$ with
a nondegenerate even pairing $V_{ij} \times V_{ji} \rightarrow k$
making $V_{ij} \cong V_{ji}^*$.  In the event $i = j$, we define a
single graded vector space $V_{i}$ together with an even nondegenerate
symmetric bilinear form $V_i \times V_i \rightarrow k$ making $V_i
\cong V_i^*$.

Then, we have products $m_{n,i_1, \ldots, i_{n+1}}: V_{i_1 i_2} \o
V_{i_2 i_3} \o \cdots \o V_{i_n i_{n+1}} \rightarrow V_{1 i_{n+1}}$ for
every choice of indices $i_1, \ldots, i_{n+1}$ such that $i, i+1$ are
adjacent in $G$ for each $i$, and so are $1$ and $n+1$.  We define
from the pairing $V_{1 i_{n+1}} \cong V_{i_{n+1} 1}^*$ the maps $\tilde
m_{n, i_1, \ldots, i_{n+1}}: V_{i_1 i_2} \otimes V_{i_2 i_3} \o \cdots
\o V_{i_{n+1} i_1} \rightarrow k$.

The identities these maps are required to satisfy are then exactly the same
as the ones for an $A_\infty$ algebra, in cases when the identities apply.

The above structure can also be viewed as a genuine
 $A_\infty$ category such that hom-groups between non-adjacent vertices
 are zero.

Given such a structure, the construction from before just gives classes
in the stable homology $PH(L_{\infty G})$, just as before:
\begin{thm} \label{csht}
Given any $A_\infty$ category $\mathcal A$ whose objects are the set
 $I$ of vertices of $G$, such that $\Hom(i,j) = 0$ if $i$ is not
 adjacent to $j$ in $G$, and with all $m_1$'s equal to zero, one can
 explicitly construct a stable homology class in $PH_k(L_{\infty G})$
 for each $k$ using Kontsevich's method.  In other words, one can
 construct a class in $H_c^{2g-2+m+k} (M^{\text{comb}}_{g,m,G,X})$ for any $k, g,$ and $m$, from any such
 $A_\infty$ category.
\end{thm}

\section{Moyal quantization of necklace Lie algebras}\label{intro2}
\subsection{Reminder on Moyal product.}\label{reminder} Let $V$ be a finite
dimensional vector space equipped with a nondegenerate bivector
$\pi\in\wedge^2V$. Associated with $\pi$ is a Poisson bracket
$f,g\mapsto \{f,g\}:=\langle df\wedge dg,\pi\rangle$ on $\k[V],$ the
polynomial algebra on $V$.  The usual commutative product $m: \k[V]\o
\k[V]\to\k[V]$ and the Poisson bracket $\{-,-\}$ make $\k[V]$ a Poisson
algebra.  This Poisson algebra has a well-known {\em Moyal-Weyl
quantization} (\cite{M}, see also \cite{CP}).  This is an associative
star-product depending on a formal quantization parameter $\h$, defined
by the formula \beq{star} f *_\h g:=m\ccirc e^{\frac{1}{2} \h \pi} (f\o
g)\in \k[V][\h],\quad \forall f,g\in\k[V][\h].  \eeq

To explain the meaning of this formula, view elements of $\Sym V$
as
constant-coefficient differential operators on $V$. Hence,
an element of $\Sym V\o \Sym V$ acts  as a
constant-coefficient differential operator on 
the algebra $\k[V]\o\k[V]=\k[V\times V].$ Now, identify
$\wedge^2V$ with the subspace of skew-symmetric tensors in $V\o V$.
This way, the bivector $\pi\in\wedge^2V\subset V\o V$ becomes
a second order constant-coefficient differential operator
$\pi: \k[V]\o\k[V]\to\k[V]\o\k[V].$ Further, it is clear that
for any element $f\o g\in\k[V]\o\k[V]$ of total degree
$\leq N$, all terms with $d>N$ in the
infinite sum $e^{\frac{1}{2} \h\cdot \pi}(f\o g)=\sum_{d=0}^\infty
\frac{\h^d}{2^d d!} \pi^d(f\o g)$
vanish, so the sum makes sense.
Thus,
the 
symbol $m\ccirc e^{\frac{1}{2} \h\cdot \pi}$ in the right-hand side
of formula \eqref{star} stands for the composition
$$ \k[V]\o\k[V]\stackrel{e^{\frac{1}{2} \h\cdot \pi}}\tooo \k[V]\o\k[V]\o\k[\h]
\stackrel{m\o\Id_{\k[\h]}}\tooo  \k[V]\o\k[\h],
$$
where
$e^{\frac{1}{2} \h\cdot \pi}$
is an infinite-order formal differential operator.

In down-to-earth terms, choose coordinates
$x_1, \ldots, x_n, y_1, \ldots, y_n$ on $V$ such that
the bivector $\pi$, resp., the Poisson bracket $\{-,-\}$, takes 
the canonical form
\beq{pois}
\pi = \sum_i \frac{\partial}{\partial x_i} \o \frac{\partial}{\partial
y_i} - \frac{\partial}{\partial y_i} \o \frac{\partial}{\partial x_i},
\quad\text{resp.,}\quad
\{f,g\}=\sum_i \frac{\partial f}{\partial x_i}\frac{\partial g}{\partial
y_i} - \frac{\partial  f}{\partial y_i} \frac{\partial g}{\partial x_i}.
\eeq
Thus, in canonical coordinates
$x=(x_1, \ldots, x_n), y=(y_1, \ldots, y_n),$ formula \eqref{star} for the Moyal
product
reads
\begin{align}\label{star1}
(f *_\h g)(x,y)&=\sum_{d=0}^\infty
\frac{\h^d}{2^d d!}\left(
\sum_i \frac{\partial}{\partial
x'_i} \frac{\partial}{\partial
y''_i} - \frac{\partial}{\partial y'_i} 
\frac{\partial}{\partial x''_i}
\right)^df(x',y') g(x'',y'')\Big|_{{x'=x''=x}\atop{y'=y''=y}}\nonumber\\
&=\sum_{\mathbf{j},\mathbf{l}\in\Z^n_{\geq 0}}
(-1)^{\mathbf{l}|}\frac{\h^{|\mathbf{j}|+|\mathbf{l}|}}{2^{|\mathbf{j}| + |\mathbf{l}|}
\mathbf{j}!\,\mathbf{l}!}\cdot
\frac{\partial^{\mathbf{j}+\mathbf{l}}f(x,y)}{\partial x^\mathbf{j}\partial y^\mathbf{l}}
\cdot
\frac{\partial^{\mathbf{j}+\mathbf{l}}g(x,y)}{\partial y^\mathbf{j}\partial x^\mathbf{l}},
\end{align}
where for $\mathbf{j}=(j_1, \ldots,j_n)\in \Z^n_{\geq 0}$
we put $|\mathbf{j}|=\sum_i j_i$ and given
$\mathbf{j},\mathbf{l}\in \Z^n_{\geq 0},$ write
$$\frac{1}{\mathbf{j}!\,\mathbf{l}!}\frac{\partial^{\mathbf{j}+\mathbf{l}}}
{\partial x^{\mathbf{j}}\partial y^{\mathbf{l}}}:=
\frac{1}{j_1!\ldots
j_n!l_1!\ldots l_n!}\cdot\frac{\partial^{|\mathbf{j}|+|\mathbf{l}|}}
{\partial x_1^{j_1}\ldots\partial x_n^{j_n}\partial
y_1^{l_1}\ldots\partial y_n^{l_n}}.
$$

A more conceptual approach to formulas \eqref{star}--\eqref{star1}
is obtained by introducing the {\em Weyl algebra} $A_\h(V)$.
This is a $\k[\h]$-algebra defined by the quotient
$$
A_\h(V):= (TV^*)[\h]/I(u\o u' - u'\o u-\h\cdot\langle \pi,
u\o u'\rangle)_{u,u'\in V^*},
$$
where $TV^*$ denotes the tensor algebra of the vector space $V^*$,
and $I(\ldots)$ denotes the two-sided ideal generated by the
indicated set.
Now, a version of the Poincar\'e-Birkhoff-Witt (PBW) theorem
says that the natural {\em symmetrization map}
yields a  $\k[\h]$-linear bijection
$\phi_W: \k[V][\h]\iso A_\h(V)$.
Thus, transporting the multiplication map in the Weyl algebra $A_\h(V)$ via
this bijection, one obtains an associative product
$$\k[V][\h]\o_{\k[\h]}\k[V][\h]\to\k[V][\h],
\quad f\o g\mapsto \phi_W^{-1}(\phi_W(f)\cdot \phi_W(g)).
$$
It is known that this  associative product is equal to the one
given by formulas  \eqref{star}--\eqref{star1}.

\subsection{The quiver analogue.} The second goal of this paper 
is to extend the constructions outlined above to noncommutative
symplectic geometry. Specifically, it turns out that the necklace Lie
algebra defined earlier can be expressed in a form that is analogous to
the Poisson algebra on commutative polynomials (which after all is just
the Lie algebra of derivations in the commutative world).  Then, we will
produce a quantization of the symmetric algebra of the necklace Lie
algebra given by explicit formulas analogous to formulas
\eqref{star}--\eqref{star1}.

In more detail, fix a quiver with vertex set $I$ and edge set $Q,$ and
let $\dq$ be the double of $Q$ obtained by adding
reverse edge $e^*\in\dq$ for each edge $e\in Q$.
Let $P$ be the {\em path algebra} of $\dq$. The commutator
quotient space $P/[P,P]$ may be identified naturally
with the space $L$ spanned by cyclic paths (forgetting which was
the initial edge), sometimes called
{\em necklaces}. Letting $\text{pr}_L: P \rightarrow P/[P,P] = L$ be the
projection, there is a natural bilinear
pairing 
\beq{pair}\{-,-\}:\
L\o L\to L,\quad
f\o g\mapsto\{f,g\}:= \text{pr}_L \biggl( \sum_{e \in Q} \frac{\partial f}{\partial e} 
\frac{\partial g}{\partial e^*} - \frac{\partial f}{\partial e^*} 
\frac{\partial g}{\partial e} \biggr).
\eeq
For this to make sense, we interpret
$\frac{\partial}{\partial e}, \frac{\partial}{\partial e^*}$
appropriately as maps $L \rightarrow P, P \rightarrow P$, using the formula
$\frac{\partial}{\partial e} (a_1 \ldots a_n) = \sum_{a_r = e} a_{r+1}
a_{r+2} \cdots a_n a_1 \cdots a_{r-1}$.
Then, this formula is a quiver analogue of \eqref{pois}, and provides $L$
with a Lie algebra structure identical with the necklace Lie
bracket defined earlier (this is easy to check; see \cite{G}, \cite{BLB}).
More recently, the second author showed in \cite{S}
that there is also a natural Lie {\em cobracket}
on $L$.
To explain this, write $a_1\cdots a_p\in P$ for a path
of length $p$ and let
 $1_i$ denote the trivial (idempotent) path at the
vertex $i\in I$. Further, for
 any edge $e\in \dq$ with head $h(e)\in I$ and tail 
$t(e)\in I$, 
 let
$D_e: P\to P\o P$ be
the derivation defined by the assignment
$$D_e:\
P\to P\o P,\quad
a_1\cdots a_p\mapsto\sum_{a_r=e} 
a_1\cdots a_{r-1}1_{t(e)}\o1_{h(e)} a_{r+1}\cdots a_p.
$$
The map $D_e$ is a derivation. Moreover, the
following map,  cf. \cite[(1.7)-(1.8)]{S}:
\beq{delta}
\delta: L\to L\wedge L,
\quad
f\mapsto \delta(f)= (\text{pr}_L \o \text{pr}_L) 
\biggl( \sum_{e \in Q} D_e(\frac{\partial f}{\partial e^*})
- D_{e^*}(\frac{\partial f}{\partial e}) \biggr)
\eeq
(that is, in a sense, dual to \eqref{pair})
makes the Lie algebra $L$ a Lie {\em bialgebra},
to be referred to as the {\em necklace Lie bialgebra}.

The necklace Lie bialgebra admits a very interesting quantization.
Specifically, the main construction of \cite{S} produces
 a Hopf $\k[\h]$-algebra
$A_\h(Q)$ equipped with an algebra
isomorphism $A_\h(Q)/\h\cdot A_\h(Q)\iso \Sym L,\,f\mapsto\pr{f}.$
The algebra $A_\h(Q)$ is a quantization of the
Lie bialgebra $L$ in the sense that $A_\h(Q)$ is flat over $\k[\h]$ and, 
for any $a,b\in A_\h(Q),$ one has
$$
\pr\left(\frac{ab-ba}{\h}\right)=\{\pr{a},\pr{b}\},\quad
\text{and}\quad
\pr\left(\frac{\Delta(a)-\Delta^{op}(a)}{\h}\right)=\delta(\pr(a)),
$$
where $\Delta: A_\h(Q)\to A_\h(Q)\o_{\k[\h]}A_\h(Q)$ denotes the
coproduct in the Hopf algebra $A_\h(Q)$,
and where $\Delta^{op}$ stands for the map
$\Delta$ composed with the flip of the two factors in
$A_\h(Q)\o_{\k[\h]}A_\h(Q).$

\subsection{Moyal quantization for quivers.} In  \cite{S}, the Hopf algebra
$A_\h(Q)$ was defined, essentially, by generators
and relations. Thus, the algebra
$A_\h(Q)$  may be viewed, roughly, as a quiver analog 
of the Weyl algebra $A_\h(V)$. One of the main
results proved in  \cite{S} is a version of
the PBW theorem. The PBW theorem
 insures 
that $A_\h(Q)$ is isomorphic
to $(\Sym L)[\h]$ as a $\k[\h]$-module, in particular, it is
flat over $\k[\h]$. 

One goal of the present paper is to provide an alternative construction
of the Hopf algebra $A_\h(Q)$. Instead of defining the algebra by
generators and relations, we define a multiplication $m$ and
comultiplication $\Delta$ on the vector space $(\Sym L)[\h]$ by explicit
formulas which are both analogous to formula \eqref{star} for the Moyal
star-product. In fact, suitably interpreted, they will be written as $f
*_\h g = e^{\frac{1}{2} \h \pi}(f \o g)$ and $\Delta_h(f) =
e^{\frac{1}{2} \h \pi'} f$.  We explain how to directly check
associativity, coassociativity and compatibility of $m$ and $\Delta$,
yielding an approach (up to difficulties involving the antipode)
independent of that used in \cite{S}.

Further, in complete analogy with the case of Moyal-Weyl quantization,
we construct a symmetrization map $\Phi: (\Sym L)[\h]\to A_\h(Q)$. This
map is a bijection, and we show that Hopf algebra structure on $(\Sym
L)[\h]$ defined in this paper may be obtained by transporting the Hopf
algebra structure on $A_\h(Q)$ defined in \cite{S} via $\Phi$.

\subsection{Representations for the Moyal quantization.}
In \cite{G}, an interesting representation of the necklace Lie algebra is presented
which is quantized in \cite{S}.  Namely, for any representation of the
double quiver $\overline{Q}$ assigning to each arrow $e \in \overline{Q}$ the
matrix $M_e: V_{t(e)} \rightarrow V_{h(e)}$, we can consider the map $L \rightarrow
\k$ given by $e_1 e_2 \cdots e_m \mapsto \tr(M_{e_1} M_{e_2} \cdots M_{e_m})$.  More
generally, if $\mathbf{l} \in \Z_{\geq 0}^I$, then we can consider the representation
space $\mathrm{Rep}_{\mathbf{l}}(\dq)$ of representations with dimension vector $\mathbf l$,
meaning that $\mathrm{dim}\ V_i = l_i$.  Then this is a vector space of dimension
$\sum_{e \in \dq} l_{t(e)} l_{h(e)}$.  It has a natural bivector $\pi((M_{e})_{ij},
(M_{f})_{kl}) = \delta_{il} \delta_{jk} [e,f]$, where $[e,f] = 1$ if $e \in Q, f = e^*$
and $[e,f] = -1$ if $f \in Q, e = f^*$, with $[e,f]=0$ otherwise. We then have 
the Poisson algebra homomorphism
\beq{trrep}
\tr_{\mathbf{l}}: \Sym L \rightarrow \k[\mathrm{Rep}_{\mathbf{l}}(\dq)], \quad 
\tr_{\mathbf{l}}(e_1 e_2 \cdots e_m)(\psi) = \tr(M_{e_1} M_{e_2} \cdots M_{e_m}).
\eeq

In \cite{S}, this representation was quantized by a representation $\rho_\lb: A \rightarrow
\mathcal D(\mathrm{Rep}_\lb(Q))$, where the latter is the space of differential operators
with polynomial coefficients on $\mathrm{Rep}_\lb(Q)$.  We may modify the $\rho_\lb$ 
and $A$ slightly to obtain $\rho_\lb^\h, A_\h$ so that we have the following diagram:
\begin{equation}
\xymatrix{
\Sym L \ar[rr]_-{\mathrm{asympt. inj.}}^{\tr_\d} & & \k[\Rep_\d(\dq)] \\
A_\h \ar[u] \ar[rr]_-{\mathrm{asympt. inj.}}^{\rho_\d^\h} & &
D_Q
\ar[u]}
\end{equation}
Here, $A_\h$ is obtained from $A$ by modifying (3.3) in \cite{S} so that the right-hand
side has an $\h$ just like (3.4). [Note: More generally, it makes sense to consider
the space where (3.3) has an independent formal parameter $\hbar$; for the Moyal
version, we want the two to be the same.] Then, the representations $\rho_\lb^\h$
send elements $(e_1, 1) (e_2, 2) \cdots (e_m, m) \in A_\h$ (see \cite{S}: this is one
lift of $e_1 e_2 \cdots e_m \in L$)
to operators $\sum_{i_1, i_2, \cdots, i_m} 
\iota(e_1)_{i_1 i_2} \iota (e_2)_{i_2 i_3} \cdots \iota(e_m)_{i_m i_1}$, 
where $\iota(e)$ is the matrix $M_e$ if $e \in Q$, and $\iota(e^*) = M_{e^*}$ for
$e \in Q$, where $M_{e^*}$ is the matrix given by 
$(M_{e^*})_{ij} = -\h \frac{\partial}{\partial (M_e)_{ji}}$.  Then, the space $D_Q \subset
\mathcal D(\mathrm{Rep}_\lb(Q))$ is just generated by $e_{ij}, 
-\h \frac{\partial}{\partial e_{ji}}$.

The diagram indicates that the representations are ``asymptotically injective'' in
the sense that the kernels of the representations $\rho_\lb, \tr_\lb$ have zero
intersection, and moreover, for any finite-dimensional vector subspace $W$ of the
algebra $A$, there is a vector $\lb \in N^I$ such that for each
$\lb' \geq \lb$ (i.e.~such that $l_i' \geq l_i, \forall i$, 
we have that $W \cap \text{Ker }\tr_\lb = 0$ (and similarly for $\rho$).

By construction of the map $\Phi_W$,
the Moyal quantization fits into a diagram as follows:
\begin{equation} \label{2d}
\xymatrix{ \Sym L[\h]_\moy \ar[rr]_{\mathrm{asympt. inj.}}^{\tr_\lb[\h]} \ar[d]
\ar@/_5pc/[dd]^{\Phi_W}_{\sim} & & \k[\h][\Rep_\lb(\dq)]_\moy \ar[d]
\ar@/^5pc/[dd]^{\phi_W}_{\sim} \\ \Sym L
\ar[rr]_-{\mathrm{asympt. inj.}}^{\tr_\lb} & & \k[\Rep_\lb(\dq)] \\ A_\h \ar[u]
\ar[rr]_-{\mathrm{asympt. inj.}}^{\rho^h_\lb} & & D_Q \ar[u] }
\end{equation}
Here, we denote by $\k[\h][\Rep_\lb(\dq)]_\moy$ the Moyal quantization
of $\k[\Rep_\lb(\dq)]$ using the bivector $\pi$, and by $\Sym
L[\h]_\moy$ the quiver version to be defined in this article.  Because
of the asymptotic injectivity, to prove that a Moyal quantization
exists completing the diagram, all that is necessary is the map
$\Phi_W$; then the definitions of the product, coproduct, and antipode
follow.  However, the definitions are interesting in their own right.

\subsection{Organization of the article.}
The article is organized as follows: In Section \ref{shs}, we prove
Theorem \ref{sht}, classifying the stable homology of the necklace Lie
algebra associated to any quiver. In Section \ref{cshs}, we prove
Theorem \ref{csht}, showing that classes of this stable homology can be
obtained from certain $A_\infty$ categories. 

In Section \ref{mps}, we will define the Moyal product $*_\h$ on $\Sym
L[\h]$. In Section \ref{phiws}, we define the map $\Phi_W$.  Next, in
Section \ref{mpts}, we show that this transports the product on $A_\h$
to the product $*_\h$. Finally, in Section \ref{ass}, we directly prove
the associativity of $*_\h$.

In Section \ref{cps} we define the Moyal coproduct $\Delta_\h$.  Then,
in Section \ref{cpts}, we show that $\Delta_\h$ is obtained by
transporting the coproduct from $A_\h$ using $\Phi_W$.  Section
\ref{casss} proves directly that $\Delta_\h$ is coassociative.

In Section \ref{as} we give the definition of antipode $S$, which
clearly is the one obtained from $A_\h$ by transportation.  This makes
$\Sym L[h]_\moy$ a Hopf algebra satisfying $S^2 = \Id$.  The
eigenvectors of $S$ are just products of necklaces, with eigenvalue $\pm
1$ depending on the parity of the number of necklaces.

Note that we have a direct proof (see \cite{GS}) of the compatibility of
the product and coproduct (the bialgebra condition), but have omitted it
to save space (since the compatibility follows immediately from \cite{S}
using the comparison).  We do not know of a direct proof that the
formula we give for antipode in Section \ref{as} satisfies the antipode
condition (although it follows from \cite{S}).

We will make use of the following tensor convention throughout:
\begin{ntn}
If $S, T$ are $\k[\h]$-modules, then we will always mean by $S \o T$ the
tensor product over $\k[\h]$ (never over $\k$).
\end{ntn}

\begin{rem} It should also be possible to give formulas for Euler
 characteristic of the quiver versions of stable homology and moduli
 space, following the well-known orbifold results \cite{HZ}, and the
 formula for a sum of Euler characteristics of the topological spaces
 given in \cite{GK}.  This will hopefully be the subject of an upcoming
 paper.
\end{rem}

\section{Stable Homology} \label{shs}
In this section we prove Theorem \ref{sht}.  Most of the
 material is a straightforward generalization of the proof in \cite{K},
 whose details are spelled out in references such as \cite{CV}. Thus, we
 endeavor to be as concise as possible.

\subsection{Reminder on the ribbon graph complex.} \label{rgcss}
  A ribbon graph is a triple $\Gamma = (H, \iota,
\gamma)$ where $H$ is a set of ``half-edges'', $\iota: H \rightarrow H$
is a fixed-point-free involution, and $\gamma: H \rightarrow H$ is a
permutation.  The orbits of $\iota$ are called the edges $E$, and the
orbits of $\gamma$ are called the vertices $V$ of $H$.  Also, we will
call the orbits of $\gamma \circ \iota$ the faces $F$ of the graph. Visually,
$\iota$ interchanges two halves of each edge (this can also be viewed as
changing the orientation of the edge), whereas $\gamma$ cyclically
rotates the half-edges at each vertex.  

\begin{ntn} For each vertex $v \in V$, let $H_v \subset H$ denote the
 set of half-edges making up the corresponding $\gamma$-orbit.  For each
 edge $e \in E$, let $H_e \subset H$ denote the set of half edges making
 up the corresponding $\iota$-orbit (of size two).
\end{ntn}

In other words, a ribbon graph is a undirected graph endowed with a
cyclic ordering of the half-edges that meet at each vertex (so, one
half-edge for each edge whose endpoints are that vertex and another one,
and two half-edges for each loop at that vertex).  

This is called a ribbon (or fat) graph because it can also be viewed as
a graph whose edges have a thickness, so that the cyclic ordering of
half-edges at each vertex is the condition that a suitable neighborhood
of the vertex must be homeomorphic to a thick star.  It is clear
that the definition of ribbon graph given above is equivalent to
homeomorphism classes of such thickened graphs with labeled vertices and
edges, so that vertices and edges map to each other.

We will restrict our attention to \textbf{connected} ribbon graphs,
which means that the thickened graph is a connected topological space,
or that $\gamma$ and $\iota$ together act transitively on $H$. From now
on, we will assume that our graph is connected.

For a given ribbon graph $\Gamma$, let $g = g(\Gamma) = 1 - \frac{1}{2}
(\#(V) - \#(E) + \#(F))$ be the \textbf{genus} of the graph, and let $n
= n(\Gamma) = \#(F)$ equal the number of faces.  It is clear that the
thickened graph is homeomorphic to a genus-$g$ surface with $n$
punctures.

The complex $M_{g,m}^{\text{comb}}$ is constructed by adding one orbicell
$C_\Gamma := C_{[\Gamma]}$ for each isomorphism class $[\Gamma]$ of
ribbon graphs $\Gamma$ of genus $g$ with $n$ faces, such that all
vertices have valence $\geq 3$.  Picking a representative $\Gamma$ of
the
isomorphism class, the cell is a copy of $\R_+^{E} /
\text{Aut}(\Gamma)$, which can be considered as an orbifold or as a
topological space.  It corresponds to choosing the lengths of the edges,
which determines uniquely the complex structure of the thickened graph.

The dimension of such an orbicell is the number of edges.  As a consequence
of the valence condition, the maximal-dimensional cells are those where
all vertices have valence $3$, and in this case, the dimension is
$6g-6+3n = 3 \#(E) - 3 \#(V) = \#(E)$.

Finally, in order to define the differential, it is necessary to
introduce an orientation.  This can be done in several equivalent ways;
the simplest for our purposes is Kontsevich's original definition:

\begin{defn}
An orientation of a ribbon graph is an orientation on the real vector
 space $\R^E \oplus \R^F$. Here $\R^X$ is defined to be the real vector
 space with basis $X$.  An orientation on a vector space $W$ is just an
 element of $((\det W) \setminus 0) / \R_+$, i.e.~a choice of
 identification $\det W \cong \R$ up to scaling by a positive factor.
\end{defn}

\begin{defn}
A ribbon graph is \textbf{orientable} if there is no automorphism
 reversing the orientation.
\end{defn}

The \textbf{ribbon graph complex}, $RG^{g,m}_\cdot$, is defined,
\cite{K}, as a vector space with a basis consisting of isomorphism
classes of certain connected ribbon graphs $\Gamma$ with orientation $\ort$,
modded by the relation $(\Gamma, \ort) = - (\Gamma, -\ort)$ (in
particular this kills ``nonorientable'' graphs).  The allowed ribbon
graphs are those whose vertices have valence $\geq 3$, which have $m$
faces and genus $g$.  The differential is given by $d(\Gamma, \ort) =
\sum_{e \in E(\Gamma)} (\Gamma/e, \ort_e)$, where $\ort_e$ is
determined from $e$ by sending $\omega_1 \otimes \omega_2 \in
\text{det}(\R^E) \otimes \text{det}(\R^F) = \det (\R^E \oplus \R^F)$
to $e^*(\omega_1) \otimes \omega_2$, where $e^*$ is the element of the
dual basis to $E$ corresponding to $e$.

Finally, the degree is given by the number of edges: so $RG^{g,m}_\ell$
is spanned by ribbon graphs of genus $g$ with $m$ faces and $\ell$ edges.

\begin{rem}
This construction generalizes to graph complexes generated by any cyclic
operad \cite{CV}, replacing the set of faces by $H^1(\Gamma, \R)$.
\end{rem}

The space $M_{g,m}^{\text{comb}}$ is glued together such that the
chain complex of the ribbon graph complex, with trivial coefficients
in a characteristic zero field $k$, is isomorphic to a
chain complex associated with the orbicell complex
$M_{g,m}^{\text{comb}}$.  (Note that this
requires that the nonorientable graphs vanish; see Remark \ref{orr}.)
Then, the well-known theorem \cite{P}, \cite{H} says that
$M_{g,m}^{\text{comb}} \iso M_{g,m} \times \R^{n}$ by thickening graphs
(and mapping to $\R^n$ the perimeters of the labeled faces).

To define this precisely, consider an ordering of the vertices $v_1,
\ldots, v_{\#(V)}$. This determines an element $M = \tilde m_{\#(v_1)}
\otimes \cdots \otimes \tilde m_{\#(v_\#(V))}$.  Furthermore, consider a
choice of ordering of the half-edges of each vertex that is in cyclic
order ($(h, \gamma(h), \gamma^2(h), \ldots)$), i.e.~a choice of initial
half-edge. This is called a choice of \textbf{ciliation} of each
vertex.  This assigns to the $2\#(E)$ components of $M$ a fixed labeling
by the half-edges $H$.  Next, consider a choice of ordering of the
edges, and a choice of orientation of each edge (a choice of half-edge
for each edge).  This assigns to the $2\#(E)$ components of $C^{\otimes
\#(E)}$ a fixed labeling by the half-edges $H$.  Then, one performs the
graded contraction $\langle M, C^{\otimes \#(E)} \rangle$ which assigns
to each half-edge on the left the corresponding half-edge on the right.
This defines the weight $W(\Gamma, \ort)$, up to sign.

For a given choice of orientation $\ort$ of $\Gamma$, one can define a
natural sign choice of $W(\Gamma, \ort)$ such that $W(\Gamma, \ort) = -
W(\Gamma, -\ort)$, and such that the resulting weight is independent of
the choices made above. This follows from the following reformulation of
orientation:

\begin{prop}\cite{CV} The orientation of a graph is naturally identified
 with a choice of orientation of $\bigotimes_{v \in V} \det (\R^{H_v})
 \oplus \R^{V_e}$, where $V_e$ is the set of vertices of even
 valence.
\end{prop}

By the proposition, an orientation for a graph is naturally identified
with a choice of ordering of the half-edges of each vertex and ordering
of the set of vertices of odd valence, modding by even permutations.  
But when we restrict the orderings of the half-edges in a
vertex to ciliations, we get a canonical one for the odd-valence
vertices, so an orientation simply gives a choice of ciliation of the
even-valence vertices, modding by applying $\gamma+1$ to the half-edges
which make up any given even-valence vertex (that is, in the quotient,
applying $\gamma$ to the half-edges of a vertex is the same as flipping
the orientation).

This gives a canonical choice of sign of $W(\gamma, \ort)$: since $C$ is
even, the choice of ordering of the edges is already irrelevant to the
graded contraction.  Also, the orientation of edges is irrelevant, since
$C$ is graded symmetric.  For each swap of vertices, one gets a minus
sign iff both vertices had even valence (by parity of $\tilde m_v$),
otherwise no sign changes.  For each cyclic rotation of the ciliation of
a vertex $v$, one obtains a minus sign iff the vertex has even valence
(i.e.~the cyclic permutation is odd).  Hence, the orientation of a graph
gives a natural choice of $W(\Gamma, \ort)$, and we get a natural element
$W(\Gamma, \ort) (\Gamma, \ort)$.

The cycle in $C_j(M_{g,m}^{\text{comb}})$ that Kontsevich defines is
then just $\sum_{\Gamma} W(\Gamma, \ort) (\Gamma, \ort)$ where the sum is
over isomorphism classes of orientable graphs of genus $g$ with $n$
holes and $j$ edges; we already saw that $W(\Gamma, \ort) (\Gamma, \ort)$
does not depend on the choice of orientation.  

It is a theorem (\cite{KFD}, \cite{PS}, \cite{Pk}) that this element is
indeed a cycle.  The proof involves just summing over the weights of
oriented graphs which expand $(\Gamma, \ort)$ by one edge, and obtaining
zero.  If one restricts attention to adding an edge at one vertex, then
the weights obtained already add to zero, which is basically a direct
consequence of \eqref{ainfid}.

We state now a third characterization of orientation that we will need
for the proof:
\begin{prop} \label{thop}
The orientation of a graph is naturally identified with
$\R^{V} \oplus \bigotimes_{e \in E}
 \det(\R^{H_e})$.
\end{prop}

\begin{rem} If we were to shift the complex so as to 
make $\tilde m_n$ graded cyclically symmetric and odd for all $n$, and
 to make $C$ graded antisymmetric, then the above all works the same
 except we use Proposition \ref{thop} to characterize orientation.
 Then it is the ordering of the vertices and the orientation of the
 edges that matter, which works for the same reason as the above.
\end{rem}

\begin{rem} \label{orr}
Note that not all ribbon graphs, with valences $\geq 3$, are oriented.
 In much of the literature this goes unmentioned. For example, the
 planar ribbon graph with two vertices and four edges is
 nonorientable. Nevertheless, with rational coefficients (which is
 needed to avoid dealing with orbifold issues), such chains are equal to
 their negative (for example, as currents). 
\end{rem}

\begin{rem}
It
is pointed out  in \cite{C} that this decomposition is not necessarily natural (there
are various triangulations), so it would be nice to get this description
 without explicit use of ribbon graphs (and the Taylor coefficients $m_n$).
\end{rem}

\begin{rem} 
Note that it is actually more natural to view the result as an
isomorphism of stable \textbf{co}homology with ribbon graph cohomology.
We will see later that the result is proved by giving an explicit,
natural isomorphism of cochain complexes.  Since we work over the
rationals, this also gives the desired result for stable homology, and
we state it in this form since it is the way it has been stated since
\cite{K}.
\end{rem}

\subsection{Stable Lie algebra homology.}
Let $Q$ be any quiver.  As we recall, the ``stable quiver'' $\infty Q$
and the stable Lie algebra $L_{\infty Q}$ only depends on the undirected
graph $G$
obtained from $Q$ by reducing all edge multiplicities to $0$ or $1$.
As in \cite{K}, \cite{CV}, we compute the stable homology
$PH_k(L_{\infty G})$ by using the standard complex 
$C_n(\mathfrak g) = \Lambda^n\mfg$, with $d_n: C_n \rightarrow C_{n-1}$
given by $d_n(g_1 \wedge \cdots \wedge g_n) = \sum_{i < j} (-1)^{i+j+1} 
[g_i, g_j]
\wedge g_1 \wedge \cdots \wedge \hat g_i \wedge \cdots \wedge \hat g_j
\wedge \cdots \wedge g_n$.  Now, letting $\mfg[m]$ be the $m$-th graded
piece in $\mfg$ (all graded pieces are finite-dimensional), one can apply
Kontsevich's trick \cite{K} (generalized to arbitrary cyclic 
operads in \cite{CV}) to pass to the $\text{ad }
\mfg[2]$-invariants. As mentioned in \cite{K}, it is well known that $\mfg$
acts trivially on homology by the adjoint action.  

This is where the new work has to be done in our case: to analyze the
structure of $\mfg[2]$. It is not very difficult:
\begin{defn}\label{lqigd} 
Let $L_Q[i]$ be the $i$-th graded component of $L_Q$, with
 respect to lengths of paths.  Similarly, let $P_Q[i]$ be the $i$-th
 graded piece of the path algebra $P_Q$.  (Note that $L_{Q,+} = \bigoplus_{i \geq 2} L_Q[i]$.)
\end{defn}
The following lemma simply recalls some obvious facts that we will use
throughout:
\begin{lemma} \cite{G}, \cite{BLB}
(a) Each component $L_Q[i]$ and $P_Q[i]$ is finite
 dimensional. (b) $P_Q[0] \cong k^I$ is a semisimple ring on primitive idempotents $I$, the
 set of vertices of $Q$. Also $L_Q[0] \cong k^I$ as a vector space. 
(c) $P_Q[1]$ is a vector space of dimension
 $\#(\dq)$, with basis the arrows of the double quiver $\dq$. (d) For
 any $i, j \in I$, $i P_Q[1] j$ has as a basis those edges beginning at $i$
 and ending at $j$.  (e) We have $L_Q = (T_{P_Q^0} P_Q[1]) / [T_{P_Q^0} P_Q[1],
 T_{P_Q^0} P_Q[1]]$ as a vector space.
\end{lemma}
\begin{defn} Let $E_{ij} := i P_Q[1] j$ be the vector subpace of the path
algebra with basis
those arrows from $i$ to $j$.  Let $L_{ij} := [E_{ij} E_{ji}] \subset L_Q[2]$
 be the subspace of cyclic paths given by going from $i$ to $j$ and then
 back to $i$ (equivalently switching $i$ and $j$).
\end{defn}
\begin{ntn} When $X_{ij} = X_{ji}$
 (identically) then we will define $X_{\{i,j\}} := X_{ij} = X_{ji}$, and 
as an abuse of notation, given a set $\{i,j\}$ we will let $X_{ij}$
 denote $X_{\{i, j\}}$ (even though there is no natural choice of which
 element to label $i$ and which to label $j$).
\end{ntn}
\begin{defn}\cite{G} Let $\omega$ be the natural symplectic form on
 $L_Q[1]$ given by $$\omega(e, f) = \begin{cases} 1, & \text{if } e \in
				   Q, f = e^*; \\
				   -1, & \text{if } f \in Q, e = f^*; \\
				   0, & \text{otherwise.}
\end{cases}
$$
\end{defn}
\begin{lemma} Using the isomorphisms $E_{ij} \cong E_{ji}^*$ induced by
$\omega$, each $L_{ij}$ is a sub-Lie algebra of $L_Q[2]$ isomorphic to
 $\mathfrak{gl}(E_{ij})$ and each $L_{ii}$ is a sub-Lie algebra
 isomorphic to $\mathfrak{sp}(E_{ii})$. In particular, $L_Q[2]$ is
 semisimple.
\end{lemma}
\begin{proof}
The fact that $\omega$ induces $E_{ij} \cong E_{ji}^*$ follows from the
 definition of $\omega$.  
Note that we can restate the bracket on $L_Q[2]$ as
 $[ab, cd] = \omega(a,c) bd + \omega(a, d) bc + \omega(b,c) ad +
 \omega(bd) ac$.  Using the isomorphisms $E_{ij} \cong E_{ji}^*$, the
result easily follows (for the case $i \neq j$, two of the terms in the
 above vanish; for the case $i = j$, it is easiest to note that $L_{ii}$
 acts on $E_{ii}$ by $(ab) * (c) = \omega(a,c) b + \omega(b,c) a$.)  The
 fact that $L_Q[2]$ is semisimple follows because we now can decompose
 $L_Q[2]$, as a Lie algebra, into the sum of the Lie ideals $L_{ij}$
 which are simple.
\end{proof}

The following two lemmas are obvious so their proofs are omitted.
\begin{lemma} The space $P_Q^m$ decomposes as an $L_Q[2]$-module (under
 the adjoint action) as follows:
\begin{equation}
P_Q^m = \bigoplus_{i_1, \ldots, i_{m+1} \in I} E_{i_1 i_2} E_{i_2 i_3}
\cdots E_{i_m i_{m+1}}.
\end{equation}
Similarly, the space $L_Q^m$  decomposes as an $L_Q[2]$-module as
\begin{equation}
L_Q^m = \bigoplus_{(i_1, \ldots, i_m) \in I^{m} / \Z/m} 
[E_{i_1 i_2} E_{i_2 i_3} \cdots E_{i_m i_{1}}].
\end{equation}
Here $I^m/\Z/m$ denotes the cyclic $m$-tuples of elements of $I$, with
 no preferred starting element.
\end{lemma}
\begin{defn} Let $CP_Q^m = \sum_{i \in I} i P_Q^m i \subset P_Q^m$
 denote the space of \textbf{closed} paths.
\end{defn}
\begin{lemma} \label{cpll}
The group $\Z/m \subset \Sigma_m$ acts canonically on 
$CP_Q^m$  (by cyclic permutations) with quotient $L_Q^m$. We have
\begin{equation}
P_Q^m = \bigoplus_{i_1, \ldots, i_{m} \in I} E_{i_1 i_2} E_{i_2 i_3}
\cdots E_{i_m i_{m}}.
\end{equation}
\end{lemma}
To study the cohomology of $L_Q^m$, we will need to make use of the
fact that $L_Q[2]$ is a finite-dimensional reductive Lie algebra
in characteristic zero, so that one has natural isomorphisms
$V_{L_Q[2]} \cong V^{L_Q[2]}$ and $(V^*)_{L_Q[2]} \cong (V_{L_Q[2]})^*$.
The same holds replacing $L_Q[2]$ by finite groups.
\begin{lemma} \label{cpi}
A basis for the invariants $((CP_Q^m)^*)^{L_Q[2]}$ is labeled by the
 pairs $(T, \beta),$ where $T$ is a fixed-point-free involution of 
$\{1, 2, \ldots, m\}$, and $\beta: \{1, 2, \ldots, m\} \rightarrow I
 \times I$ is a 
 map such that:  (i) $\beta \circ T = \sigma \circ \beta$, 
 where $\sigma(i,j) = (j,i)$ is the flip; (ii) $\beta(a)_2 =
 \beta(a+1)_1$ for any $a \in \Z/m$ (using addition in $\Z/m$),
where the subscript indicates the component of the pair in $I \times I$.

The corresponding basis element $v_{T, \beta}$
 is defined to be 
\begin{equation} \label{oe}
v_{T, \beta} = \prod_{a < T(a)} \omega^{a T(a)} \biggl|_{E_{a T(a)}
 \oplus E_{T(a) a}},
\end{equation}
where $\omega^{a T(a)}$ indicates that $\omega$ is applied to components
 $a$ and $T(a)$, in that order.
\end{lemma}
\begin{proof}
This follows immediately from the First Fundamental Theorem of Invariant
 Theory \cite{W}, noting that $P_Q^m = T^m_{P_Q^0} P_Q[1]$ ($m$-th tensor
 power over $P_Q^0$), and $P_Q[1] = \bigoplus_{\{i,j\} \subset I}
 E_{\{i,j\}} = E_{ij} \oplus E_{ji}$. Here
 $E_{\{i,j\}}$ is (canonically isomorphic to) the
 standard representation of $L_{ii} \cong \mathfrak{sp}(E_{ii})$ in the
 case $i =j$, and to the sum of the standard and dual standard 
representations of $L_{ij} \cong \mathfrak{gl}(E_{ij})$ in the case $i
 \neq j$.
\end{proof}
\begin{cor} \label{lqmi}
The invariants $((L_Q^m)^*)^{L_Q[2]}$
are spanned by elements 
labeled by $\Z/m$-orbits $[(T, \beta)]$ of pairs $(T,
 \beta)$ defined above.  Here, the $\Z/m$ action is the one induced by
 the action by translation on $\Z/m$.  The basis element $\bar
 v_{[T,\beta]}$ is just the image of $v_{T, \beta}$ under
 the quotient map $P_Q^m \onto P_Q^m/([P_Q, P_Q] \cap P_Q^m)$, for any
 element $(T, \beta)$ of the orbit $[(T, \beta)]$.  The nonzero elements
 form a basis, and an element is nonzero in the case that the orbit of
an element $v_{T, \beta}$ does not include $-v_{T, \beta}$.
\end{cor}
\begin{proof}
Note that the $L_Q[2]$ and $\Z/m$-actions on $CP_Q^m$ commute. 
We see that $((L_Q^m)^*)^{L_Q[2]} \cong (((CP_Q^m)_{\Z/m})^*)^{L_Q[2]}
 \cong ((CP_Q^m)^{L_Q[2], \Z/m})^* \cong (((CP_Q^m)^*)^{L_Q[2]})^{\Z/m}$,
 giving the first result.  For the last sentence, it is clear that the
 elements that are nonzero are linearly independent, since they
 represent linearly independent elements in $((CP_Q^m)^*)^{L_Q[2]}$.
 The condition for being nonzero follows immediately from Lemma \ref{cpi}.
%
\end{proof}
\begin{lemma} \label{scpvy}
Consider the subgroup $S_{\ell_1, r_1, \ldots, \ell_m, r_m} =
\prod_{i=1}^m \Sigma_{r_i} \times (\Z/{\ell_m})^{r_i} \subset
\Sigma_{r_1 \ell_1} \times \cdots \times \Sigma_{r_m \ell_m}$, and the
representation $V_{r_1, \ell_1, \ldots, r_m, \ell_m} :=
\boxtimes_{i=1}^m k_{\Z/\ell_i} \boxtimes \epsilon_{r_i}$, where
$k_{\Z/\ell_i}$ is the trivial representation of $\Z/\ell_i$ and
$\epsilon_{r_i}$ is the sign representation of $\Sigma_{r_i}$.  Let
$CP_{\ell_1, r_1, \ldots, \ell_m, r_m} := CP_Q^{\ell_1})^{\otimes_k r_1}
\otimes_k \cdots \otimes_k (CP_Q^{\ell_m})^{\otimes_k r_m}$.  Then
tensoring gives an isomorphism $\mathrm{Hom}(V_{\ell_1, r_1, \ldots, \ell_m,
r_m}, CP_{\ell_1, r_1, \ldots, \ell_m, r_m})$ $\iso V_{\ell_1, r_1,
\ldots, \ell_m, r_m} \otimes_{k S_{\ell_1, r_1, \ldots, \ell_m, r_m}}
CP_{\ell_1, r_1, \ldots, \ell_m, r_m}) = L_{\ell_1, r_1, \ldots, \ell_m,
r_m} := \Lambda^{r_1} L_Q^{\ell_1} \otimes_k \cdots \otimes_k
\Lambda^{r_m} L_Q^{\ell_m}$. This should be thought of as identifying
the space $L_{\ell_1, r_1, \ldots, \ell_m, r_m}$ appearing in the
standard complex for $L_Q$ with the coinvariants/invariants of
$CP_{\ell_1, r_1, \ldots, \ell_m, r_m}$ with respect to the twisted
action of $S_{\ell_1, r_1, \ldots, \ell_m, r_m}$.
\end{lemma}
\begin{proof} This follows immediately from Lemma
 \ref{cpll}.
\end{proof}
\begin{cor}$\Lambda^{r_1}(L_Q^{\ell_1}) \otimes \cdots
 \otimes \Lambda^{r_m} (L_Q^{\ell_m})$ is spanned by orbits $[(T, \beta)]$
of a fixed-point-free involution $T$ of $B_{\ell_1, r_1, \ldots, \ell_m,
 r_m} := (\Z/\ell_1) \times \{1, 2, \ldots, r_1\} \sqcup
\cdots \sqcup (\Z/\ell_m) \times \{1, \ldots, r_m\}$ 
and a map $\beta: B_{\ell_1, r_1, \ldots, \ell_m, r_m} \rightarrow I
 \times I$
 satisfying: (i) $\beta \circ T = \sigma \circ \beta$, and (ii) $\beta(a)_2 =
 \beta(a+1)_1$, for any $a \in \Z/\ell_i$, with addition modulo
 $\ell_i$.  The orbits are under the action of $S_{\ell_1, r_1, \ldots,
 \ell_m, r_m} \subset \Sigma_{\ell_1 r_1} \times \cdots \times
 \Sigma_{\ell_m r_m}$, which acts in the obvious standard way.  Taking
 the nonzero such elements gives a basis of the invariants.
\end{cor}
\begin{proof} This follows from the previous Lemma \ref{scpvy} 
and Corollary \ref{cpi}, using the fact that the Fundamental Theorem of
 Invariant Theory applies to $CP_{\ell_1, r_1, \ldots, \ell_m, r_m}$
 just as well as it did in the case $m=1, r_1=1$. We get a formula
 similar to \eqref{oe}: a product of terms $\omega^{ab}$ acting in
 distinct components, but where the components of $\omega$ can act in
 different components ($L_Q^{\ell_i}$'s) of the exterior/tensor power as
 well as in different components of an individual component $L_Q^{\ell_i}$.
\end{proof}
\begin{cor} For sufficiently large $N \in \Z_+$, a basis for the 
$L_{NQ}[2]$-invariants of $\Lambda^{r_1}
 L_Q^{\ell_1} \otimes \cdots \otimes \Lambda^{r_m} L_Q^{\ell_m}$ is
 given by isomorphism classes of orientable ribbon graphs 
together with a labeling
 of the faces by vertices (elements of $I$), such that any two adjacent
 faces have labels which are adjacent vertices in the quiver. Here
 ``orientable'' means that an orientation of the given labeled ribbon
 graph is not equivalent to the opposite orientation (as labeled ribbon graphs).
\end{cor}
\begin{proof}
The orbits $[T, \beta]$ are just isomorphism classes of ribbon graphs
 that include a label of half-edges by $I \times I$, such that the
 ``orientation-reversing'' involution $\iota$ has the flipped label, and
 such that the cyclic permutation $\gamma$ sends a label $(a,b) \in I
 \times I$ to a label $(b,c)$ for some $c \in I$.  But this is the same
 as putting a label by $I$ between any two adjacent half edges $h,
 \gamma(h)$ at each vertex of the ribbon graph, such that two adjacent
 labels correspond to adjacent vertices of the quiver, and such that the
 label between $h$ and $\gamma(h)$ is the same as the label between
 $\iota \circ \gamma(h)$ and $\gamma(\iota \circ \gamma(h))$. The latter
 just means that we have a labelling of faces of the ribbon graph by
 $I$; then the former says that adjacent faces have adjacent labels.

Note that the above orbits include an orientation, i.e.~an equivalence
 class of choice of edge-orientations and vertex-orderings (using
 Proposition \ref{thop}).  It is clear that an element is the negative
 of its reverse orientation.  On the other hand, for large enough $N$,
 the element given above does not vanish for orientable graphs, and
 they are linearly independent. This follows because, taking the
 embedding into $\Hom_{S_{\ell_1, r_1, \ldots, \ell_m,
 r_m}}(V_{\ell_1, r_1, \ldots, \ell_m, r_m}, CP_{\ell_1, r_1, \ldots,
 \ell_m, r_m})$, the image of the orientable elements (with any chosen
 orientation) are linearly independent in the case that $N$ is large
 enough (e.g.~$N$ greater than the total number of edges = $\sum_i r_i
 \ell_i$; so that, for example, each proposed basis element above has
 a nonzero term in its coordinate expansion which doesn't appear in
 any of the others.)
\end{proof}
\begin{defn} For any $m \geq 0$ and any $X \in I^{(m)}$, call an
 $X$-labeled ribbon graph a ribbon graph with $m$ faces, together with
a labeling of the faces by
 elements of $I$, such that the labels with multiplicity are given by
 the element $X \in I^{(m)}$ (unordered $m$-tuples of elements of $I$).
 An oriented $X$-labeled ribbon graph $(\Gamma, \ort)$ is just a labeled
 ribbon graph together with an orientation defined in the usual way
 (e.g.~a choice of edge orientations and an ordering of the vertices,
 modded by sign of permutations of these).
\end{defn}
\begin{defn}
Let $RG^{g,m,G,X}_\cdot$ be the complex of $X$-labeled orientable ribbon
 graphs of genus $g$ with $m$ faces, such that each vertex has valence
 $\geq 3$, where $X \in I^{(m)}$.  That is, we take the vector space
 with basis $(\Gamma, \ort)$ where $\Gamma$ is an $X$-labeled ribbon
 graph and $\ort$ is an orientation, and mod out by isomorphisms of
 ribbon graphs and the relation $(\Gamma, -\ort) = -(\Gamma,
 \ort)$. Then the differential is given by the same formula on
 oriented labeled ribbon graphs as for the usual ribbon graph complex
 $RG^{g,m}$.
\end{defn}
For technical reasons we will need the complex defined without the
valence condition:
\begin{defn}
Let $(RG')^{g,m,G,X}_{\cdot}$ be the complex defined as above but
changing the valence restriction to $\geq 2$ (rather than $\geq 3$).
\end{defn}
The reason why we prefer $RG$ is because this is a finite complex for each $g,m$ (the trivalence condition gives a bound on the number of edges so as to not exceed the given genus). One has the following result:
\begin{lemma} \label{ssl} One has $H_\cdot(RG') \cong H_\cdot(RG) \oplus
\bigoplus_{v \sim v} PH_k(\sp(2 \infty)) \oplus
 \bigoplus_{v \sim w, v \neq w}  PH_k(\gl(\infty))$.
\end{lemma}
\begin{proof} We can identify  $\bigoplus_{v \sim v} PH_k(\sp(2 \infty)) \oplus
 \bigoplus_{v \sim w, v \neq w} PH_k(\gl(\infty))$ with the homology of
the subcomplex of the complex $RG'$ where all vertices have valence $2$.
Then, it remains to see that the subcomplex of $RG'$ of graphs
containing a vertex of valence $\geq 3$ is quasi-isomorphic to its
subcomplex $RG$ (where \textbf{all} vertices have valence $\geq 3$).
This follows by a spectral sequence argument from \cite{K} (explained in
more detail in \cite{CV}).  Essentially, all one has to do is to filter
by the number of $2$-valent vertices, and the resulting spectral
sequence collapses at the first term to the homology of $RG$ because the
complex of a segment with $m$ interior vertices (all of valence $2$)
telescopes and thus has homology concentrated in degree zero.
\end{proof}

Now, we finally get to the isomorphism of cochain complexes:
\begin{lemma} For each quiver $Q$ and $N \geq 1$, 
there is a canonical epimorphism 
$$\bigoplus_{g, m, X
 \mid 2-2g-m < 0} 
C^{\cdot+2g-2+m}((RG')^{g,m,G,X}, k)
\rightarrow
 C^{\cdot}(L_{NQ}^{L_{NQ}[2]}, k)$$ of cochain complexes (with trivial
 coefficients) which is asymptotically an isomorphism (for each degree
 $j$, the map is an injective in degrees $\leq j$ for all 
$N \geq \frac{3j}{2}$). The cochain map
is given by, for any choice of orientation of the edges, ordering of
 vertices, and ciliation of the vertices of 
 $\Gamma$, $[(\Gamma, \ort)]^* \mapsto \pm \prod_{e \in E(\Gamma)}
 \omega^e$, where the sign $\pm$ is such that $\pm \ort$ is equivalent
to the orientation given by the edge orientations and vertex ordering
 (see Proposition \ref{thop}: we ignore vertex ciliations for this), 
and $\omega^e$ means that $\omega$ acts in the
 components corresponding to the edge $e$ under the choice of vertex
 ordering and ciliation, and the edge orientation says which component
 matches the first component of $\omega$.  
\end{lemma}
\begin{proof}
Well-definition: First note that the element $[(\Gamma, \ort)]^*$ 
is the dual basis element
 of $[(\Gamma, \ort)]$, using as a basis the orientable
 $(G',\eta)$-labeled ribbon graphs with any fixed chosen orientation
 (which includes $[(\Gamma, \ort)]$). Clearly the element does not
 depend on the choices of orientation. 
Also, note that the map is obviously
 independent of vertex ciliation and depends on orientation as in Proposition
 \ref{thop}. Thus it is well-defined.

Now, the fact that the map is a cochain map is not difficult: the
 cochain differential on the left sums over ways of expanding an edge of
 $\Gamma$ (since we are in the dual to the ribbon graph complex); this
 is the same as what happens on the right since, applying the image of
 $[(\Gamma, \ort)]$ (a product of $\omega^{ab}$'s) composed with the
 chain differential to $(RG')_{\cdot+1}^{g,m,G,X}$ we get exactly the sum
 of images of $[(\Gamma', \ort')]$ obtained by expanding $\Gamma$ by
 adding an edge, and not changing labels (because the expanded edge just
 corresponds to an extra $\omega$ by the definition of the necklace Lie
 algebra bracket).

The fact that the map is an epimorphism follows immediately from the
 lemmas, as does the asymptotic injectivity as stated.

Finally, for the degree computation, note that $2-2g = \#V(\Gamma) -
 \#E(\Gamma) + \#F(\Gamma)$, so that $\#V(\Gamma) = 2-2g-\#V(G')+\#E(\Gamma)$.
\end{proof}

From the above results, it follows that the primitive stable
cohomology, with trivial coefficients, is isomorphic to the homology
of $RG'$, which itself is isomorphic to the sum of the cohomology of
$M^{\text{comb}}_{g,n,G,X}$ with compact supports and the primitive
stable homologies of $\mathfrak{gl}, \mathfrak{sp}$ appearing in Lemma
\ref{ssl}.  Since we are over characteristic zero and $L_{n Q}$ has
finite-dimensional graded components for each finite quiver $Q$, there
is no problem in dualizing, so it follows that the primitive stable
homology of $L_{\infty G}$ is the same as for cohomology.  This
completes the proof of Theorem \ref{sht}.


\subsection{Proof of Theorem \ref{csht}} \label{cshs}
Consider any quiver $Q$ and the associated
multiplicity-free graph $G_Q$ constructed in Section 1 (by forgetting
multiplicities but not loops).  Let $A$ be an $A_\infty$ category with
set of objects $I$, and such that $\text{Hom}(i, j) = 0$ if $i, j \in I$
are not adjacent.

This is clearly equivalent to specifying only hom groups between
adjacent vertices and applying the axioms of an $A_\infty$ category only
to elements which involve only multiplications of adjacent vertices.

Given such a category $A$ and any $g \geq 0, m \geq 1$ with $2g-2+m \geq
0$, and any $X \in I^{(m)}$, we can associate an $L_{\infty
G}[2]$-invariant cocycle in $C_{g,m,G,X}^j$ for any degree $1 \leq j
\leq 6g-6+3m$.  Namely, we consider the sum over all isomorphism classes
of $X$-labeled ribbon graphs $\Gamma$, the element $W(\Gamma,
\ort)(\Gamma, \ort)$ defined in Section \ref{ainfcs}, which does not
depend on choice of orientation for the same reason as given in that
section, and is a cocycle also by the same computation.  This gives an
element of the stable cohomology, which is isomorphic to the stable
homology.

This completes the proof of Theorem \ref{csht}.

\section{The Moyal product}\label{mpds}
\subsection{Definition of the Moyal product $*_\h$.}\label{mps}
To define the product $*_\h$ on $\Sym L[\h]_\moy$, we proceed by
analogy: let $\pi = \sum_{e \in Q} \frac{\partial}{\partial e} \o
\frac{\partial}{\partial e^*} - \frac{\partial}{\partial e^*} \o
\frac{\partial}{\partial e}$.  For each $n \geq 0$, we define an
operator $\pi^n: \Sym L \o \Sym L \rightarrow \Sym L$, and hence
$e^{\frac{1}{2} \h \pi}: \Sym L[\h] \o \Sym L[\h] \rightarrow \Sym
L[\h]$ as follows.  We define the action of each
\begin{equation} \label{act}
T = \frac{\partial}{\partial a_1} \frac{\partial}{\partial a_2}
\cdots \frac{\partial}{\partial a_m} \o \frac{\partial}{\partial
a_1^*} \frac{\partial}{\partial a_2^*} \cdots \frac{\partial}{\partial
a_m^*}, \quad a_i \in \dq, (e^*)^* := e;
\end{equation}
and extend by linearity.  This action is best described by considering
monomials in $\Sym L$ to be collections of closed paths in $\dq$.
Each closed path corresponds to a single cyclic monomial of $L$, so a
collection of closed paths corresponds to a symmetric product of the
corresponding cyclic monomials, giving an element of $\Sym L$.  Such
elements generate all of $\Sym L$.

Take any operator of the form \eqref{act}, and two elements $P, R \in
\Sym L$, which are symmetric products (i.e.~collections) of closed
paths.  Then the element $T$ of \eqref{act} acts on $P \o R$ by
summing over all ordered choices of distinct instances of edges $e_1,
e_2, \cdots, e_m$ in the graph of $P$ such that $e_i$ is identical
with $a_i$ as elements of $\dq$, and over all ordered choices of
distinct instances of edges $f_1, f_2, \cdots, f_m$ in the graph of
$R$ such that $f_i$ is identical with $e_i^*$ as elements of $\dq$,
and adding the following element: Delete each $e_i$ from $P$ and each
$f_i$ from $R$, and join $P$ and $R$ at each $h(e_i) = t(f_i)$ and
each $h(f_i) = t(e_i)$.  The result is some element $Z \in \Sym L$
obtained from $P \o R$, which is some new collection of closed paths (or
isolated vertices, which correspond to idempotents).
So, $T(P \o R)$ is the sum of all such elements $Z$ (some of them can
be identical, of course; we are summing over the element $Z$ we get
for each choice of instances of the given edges in $P$ and $R$).

Let us explain how to make this more precise.  We define an ``abstract
edge'' to be an occurrence of an edge in a collection of necklaces 
\begin{equation} \label{pform} P =
P_1 \& \cdots \& P_k \& V,
\end{equation} 
where each $P_i = a_{i1} a_{i2} \cdots
a_{il_i}$ is a cyclic monomial (i.e.~necklace).  Here $V \in \Sym L[0]$
can be taken to be a symmetric product of vertices (paths of length
$0$).  Then, an abstract edge of $P$ is just a choice of indices
$(i,j)$.  

To sum over cuttings and gluings of two elements $P, R$ (which are
products of necklaces), we sum over collections of pairs of abstract
edges, each containing one edge of $P$ and a reverse abstract edge from
$R$, cut all edges, and glue the endpoints.

We need also to make precise what it means to ``glue the endpoints''.
Suppose $X$ is the set of abstract edges of $P$ and $Y$ the set of
abstract edges of $R$, and $I_X \subset X, I_Y \subset Y$ are the edges
we will be cutting. The possible difficulty arises in the case that
$I_X$ contains two adjacent edges: then what does it mean to glue the
endpoint that is between two edges that are both cut out?  The answer is
to define the result of the gluing to be the collection of necklaces one
gets by starting with any edge of $X$ or $Y$ which is \textbf{not} cut 
(i.e.~not in $I_X$ or $I_Y$) and define the necklace containing that
edge to be what is obtained by continually passing to the next edge,
unless it is in $I_X \cup I_Y$, in which case one passes to the edge
which follows the reverse cut edge; if this is also cut, one iterates.

To make this description precise, let $\phi: I_X \rightarrow I_Y$ be the
bijective map pairing each edge with an opposite one (i.e.~satisfying
$\pr_Y \circ \phi = * \circ \pr_X$ where $\pr$ is the projection to
$\bar Q$ and $*$ is the edge-reversal operation). Also, extend $\phi$ to
an involution of $I := I_X \cup I_Y$ (by acting by $\phi^{-1}$ on
$I_Y$).  We would like to define ``passing to the next edge'' on $X
\sqcup Y \setminus I$. To do this, first define an auxiliary map $f: X
\sqcup Y \rightarrow X \sqcup Y$ (passing to the next edge) by: if $(i,j)
\notin I_X \cup I_Y$, then $f(i,j) = (i,j+1)$ (addition taken modulo
$l_i$) (pass to the next edge from a non-cut edge); if $(i,j) \in I_X
\cup I_Y$ then $f(i,j) = \phi(i,j) + 1$ (here $(p,q)+1 :=(p,q+1)$ passes
to the next edge). That is, for cut edges, rather than passing to the
next edge, we pass to the edge following the paired cut edge.  Then, we
can define the actual ``passing to the next edge'', $f': X \cup Y
\setminus I \rightarrow X \cup I \setminus I$ by $f'(x) = f^m(x)$, where
$m \geq 1$ is the smallest positive integer such that $f^m(x) \notin I_X
\cup I_Y$.  Finally, orbits of $f'$ are just the necklaces which result
from cutting and gluing.

Now that we have defined the action of \eqref{act}, we can extend
linearly over $\k$ to obtain the action of $\pi^n: \Sym L \o \Sym L
\rightarrow \Sym L$ for any $n$, and by linearity over $\k[\h]$, also
$e^{\frac{1}{2} \h \pi}: \Sym L[\h] \o \Sym L[\h] \rightarrow \Sym
L[\h]$. (Note that only polynomials in $\h$ are required since the
application of any differential operator of degree greater than the
total number of edges appearing in a given $P \o R$ is zero).

Now, we define $*_\h: \Sym L[\h] \o \Sym L[\h] \rightarrow \Sym L[\h]$ by
\begin{equation}
P *_\h R = e^{\frac{1}{2} \h \pi} (P \o R).
\end{equation}
This defines the necessary product which allows us to define $\Sym
L[\h]_\moy$.

We can describe this more directly as follows: again let $P, R$ be of
the form \eqref{pform} with sets of abstract edges $X, Y$, respectively,
and maps $\ppi_X: X \rightarrow \dq, \ppi_Y: Y \rightarrow \dq$. 
Then 
\begin{equation}
P *_\h R = \sum_{(I_X, I_Y, \phi)} \frac{\h^{\#(I_X)}}{2^{\#(I_X)}}
s(I_X, I_Y, \phi)
PR_{I_X, I_Y, \phi},
\end{equation}
where $(I_X, I_Y, \phi)$ is any triple of a subset $I_X \subset X, I_Y
\subset Y$ and a bijection $\phi: I_X \rightarrow I_Y$ satisfying
$\ppi_Y \circ \phi = * \circ \ppi_X$, and $PR_{I_X, I_Y, \phi}$ is the
result of cutting and gluing $P$ and $R$ along this triple as described
previously.  Here $*(e) = e^*$ is the edge-reversal involution of $Q$.
The sign $s(I_X, I_Y, \phi)$ is defined by $s(I_X, I_Y, \phi) =
(-1)^{\#(I_Y \cap \ppi_Y^{-1}(Q))}$. This follows because
$e^{\frac{1}{2} \h \pi} = \sum_{N \geq 0} \frac{\h^N}{2^N}
\frac{\pi^N}{N!}$, and each $\pi^N$ involves a sum over all cuttings and
gluings of $P$ and $R$ along $N$ edges counting each ordering and
multiplying in $-1$ for each time the $\frac{\partial}{\partial e}$
appears in the second component for $e \in Q$; dividing by $N!$ means we
don't count orderings of $I_X$ so that it is only over subsets that we
sum.

In general, elements $P, R \in \Sym L[\h]$ are linear combinations over
$\k[\h]$ of such collections of necklaces, so the element $P *_\h R$
is given by summing over each choice of necklace collections in $P$
and $R$, of the product of the coefficients of the two necklace
collections and the element described in the previous paragraph.  In
other words, we sum over all ways to take the product of terms from $P$
and $R$, not just by the usual product in $\Sym L[\h]$, but also by 
$\frac{\h^p}{2^p}$ times the ways in which we can cut out $p$ matching
edges from each term and join them together (while just multiplying
the $\k[\h]$-coefficients).

\subsection{Definition of the symmetrization map $\Phi_W$.} \label{phiws}
Now, we define $\Phi_W: \Sym L[\h] \rightarrow A_\h$.  To do this, we
need to define the notion of ``height assignments''.  Let's consider a
collection of necklaces $P$ of the form \eqref{pform}.  Let $X$ be the
set of abstract edges of $P$, say $\#(X) = N$.  Then, a \textsl{height
assignment} for $P$ is defined to be a bijection $H: X \rightarrow
\{1, 2, \ldots, N\}$.  
We have the
element $P_H \in A_\h$ obtained by assigning heights to the edges in
$X$ by $H$, that is
\begin{multline}
P_H = (a_{11}, H(1,1)) \cdots (a_{1 l_1},
H((1, l_1)) \& \cdots \\ \& (a_{k 1}, H(k, 1)) \cdots (a_{k l_k}, 
H(k, l_k)) \& v_1 \& v_2 \& \cdots \& v_q.
\end{multline}
Note that we could also think of $H$ as an element of $S_N$ with some
modifications to the formula.

The element $\Phi_W$ involves taking an average over all 
height assignments:
\begin{equation}
\Phi_W(P) = \frac{1}{N!} \sum_{H} P_H,
\end{equation}
where $H$ ranges over all height assignments.  Following is the alternative
description in terms of permutations $S_N$:
Let $\theta(i,j) = j + \sum_{p = 1}^{i-1} l_p$ so that $\theta(k,
 l_k) = N$. Then
\begin{multline}
\Phi_W(a_{11} \cdots a_{1 l_1} \& a_{21} \cdots
 a_{2 l_2} \& \cdots \& a_{k 1} \cdots a_{k l_k} \& v_1 \& v_2 \& \cdots \& v_q) \\ = \sum_{\sigma
 \in S_N} \frac{1}{N!} (a_{11}, \sigma(\theta(1,1))) \cdots (a_{1 l_1},
\sigma(\theta(1, l_1))) \& \cdots \\ \& (a_{k 1}, \sigma(\theta(k, 1))) \cdots
(a_{k l_k}, \sigma(\theta(k, l_k))) \& v_1 \& v_2 \& \cdots \& v_q.
\end{multline}

\subsection{Proof that $*_\h$ is obtained from $\Phi_W$.}\label{mpts}
Let's show that $\Phi_W$ makes the diagram \eqref{2d} commute. We know
that $\Phi_W$ is an isomorphism of free $\k[\h]$-modules (using PBW
for $A_\h$, or the fact that $\rho_\d$ is asymptotically injective and
the fact that the Weyl symmetrization map is an isomorphism on the
right-hand side of \eqref{2d}).  So, once we show commutativity of the
diagram, it will follow that $\Phi_W$ induces some multiplicative
structure on $\Sym L[\h]_\moy$ making the $\Phi_W$ an isomorphism of
$\k[\h]$-algebras.  We will then want to show that this structure is
the one we have just defined, i.e.~to show that $\Phi_W$ is a
homomorphism of rings using our $*_\h$ structure.

We need to show that $\rho_\d \circ \Phi_W = \phi_W \circ \tr$. This
follows immediately from the definitions, because $\rho_\d \circ
\Phi_W$ involves summing over the symmetrization of polynomials in
$(M_e)_{ij}, \frac{\partial}{\partial (M_e)_{ji}}, e \in Q$ where
$(M_e)_{ij}$ are the coordinate functions of the matrix corresponding
to the vertex $e$; also, $\tr$ takes an element of $\Sym L[\h]_\moy$
and gives the element of $\k[\h][Rep_\d(\dq)]$ corresponding to the
trace of the (cyclic noncommutative) polynomial, which upon
substituting $(M_{e^*})_{ij} \mapsto -h \frac{\partial}{\partial (M_e)_{ji}}$
and symmetrizing (which we needed to do for this to be well-defined,
since the $(M_{e^*})_{ij}, (M_e)_{ij}$ commuted), gives the same element.

Next, let us show that the ring structure obtained from $\Phi_W$, making
$\Phi_W$ an isomorphism of rings, is exactly the product $*_\h$ we
have described in detail.
\begin{equation}
\label{pqe}
\Phi_W(P *_\h R) = \Phi_W(P) \Phi_W(R).
\end{equation}

Now we prove \eqref{pqe}.  Let's take $P = P_1 \& P_2 \& \cdots \&
P_n$, as before, to be a collection of necklaces, and similarly for $R
= R_1 \& R_2 \& \cdots \& R_m$.  (We can forget about the idempotents
such as in \eqref{pform}, since they won't affect what we have to
prove.)  Let $X$ be the set of abstract edges of $P$ and $Y$ the set
of abstract edges of $R$.  We will use $H_P: X \rightarrow \{1,
\ldots, |X|\}$ to denote a height assignment for $P$ and $H_R: Y
\rightarrow \{1, \ldots, |Y|\}$ to denote a height assignment for $R$.
Let $PR := P \& R$ denote the symmetric product of $P$ and $R$ (NOT $*_\h$.)
Let us say that a height assignment $H_{PR}: X \sqcup Y \rightarrow
\{1, \ldots, |X|+|Y|\}$ \textsl{extends} height assignments $H_P, H_R$
if $H_{PR}$ restricted to $P$ is equivalent to $H_P$ and $H_{PR}$
restricted to $R$ is equivalent to $H_R$. In other words, $H_{PR}(x_1)
< H_{PR}(x_2)$ iff $H_P(x_1) < H_P(x_2)$ for all $x_1, x_2 \in X$, and
similarly $H_{PR}(y_1) < H_{PR}(y_2)$ iff $H_R(y_1) < H_R(y_2)$ for
all $y_1, y_2 \in Y$.

Now, we know that
\begin{equation}
\Phi_W(P *_\h R) - \Phi_W(PR) = \sum_{N = 1}^\infty \frac{\h^N}{2^N} 
\Phi_W(\frac{\pi^N}{N!} (P \o R)),
\end{equation} 
and also that
\begin{multline}
\Phi_W(P) \Phi_W(R) - \Phi_W(PR) \\ = \frac{1}{(|X|+|Y|)!} 
\sum_{H_P, H_R} \sum_{H_{PR} \text{\ extending\ }H_P, H_R} (P_{H_P} *_\h 
R_{H_R}
- PR_{H_{PR}}).
\end{multline}

We are left to show, using the relations which define $A_\h$, that
\begin{equation} \label{ltsiso}
\sum_{N = 1}^\infty \frac{\h^N}{2^N} \Phi_W(\frac{\pi^N}{N!} (P \o R)) = 
\sum_{H_{PR} \text{\ extending\ }H_P, H_R} (P_{H_P} *_\h R_{H_R}
- PR_{H_{PR}})
\end{equation}

To prove this, let us fix a particular $H_P, H_R$, and $H_{PR}$, and 
expand $P_{H_P} *_\h R_{H_R} - PR_{H_{PR}}$ using
the relations that define $A_\h$.  We do this by expressing this as a sum
of commutators obtained by commuting a single edge coming from $R$
with a single edge coming from $P$.  We get
\begin{equation}
P_{H_P} R_{H_R} - PR_{H_{PR}} = \sum_{\underset{\ppi_X(x) = \ppi_Y(y)^*}{x \in X, y \in Y \text{\ such that\ 
} H_P(x) > H_R(y),} } [\ppi_X(x), \ppi_Y(y)] \h PR'_{x,y},
\end{equation}
where $PR'_{x,y}$ corresponds to taking $PR$, deleting $x$ and $y$ and
joining the endpoints, and using the unique height assignment which
gives the same total ordering on $X \sqcup Y \setminus \{x, y\}$ as the
function $H'$ given by $H'(x') = x'$ for $x' \in X$ such that $H_P(x') <
H_P(x)$, and $H'(z) = H_P(x) + H_{PR}(z)$ for all other $z \in X \sqcup
Y \setminus \{x,y\}$.  Note here that $[e,e^*] = 1$ if $e \in Q$ and
$-1$ if $e^* \in Q$.

By applying the relations repeatedly we get that
\begin{multline} \label{hpre}
P_{H_P} R_{H_R} - PR_{H_{PR}} \\ = \sum_{\underset{\text{\ such that\ }
 H_P(x_i) > H_R(y_i), \ppi_X(x_i)
= \ppi_Y(y_i)^*}{x_1, \ldots, x_k \in X, y_1,
\ldots, y_k \in Y } } [\ppi_X(x_1), \ppi_Y(y_1)] \\ \cdots [\ppi_X(x_k),
\ppi_Y(y_k)] \h^k PR''_{(x_1, y_1), \ldots, (x_k, y_k)},
\end{multline}
where $P/R''_{(x_1, y_1), \ldots, (x_k, y_k)}$ involves taking $PR$ and
deleting the pairs of edges and gluing at their respective endpoints;
and this time assigning heights by restricting $H_{PR}$ to $X \sqcup Y
\setminus \{x_1, \ldots, x_k, y_1, \ldots, y_k\}$ (and changing to a
function which has image $\{1, \ldots,
|X|+|Y|-2k\}$ which gives the same total ordering of $X \sqcup Y
\setminus \{x_1, \ldots, x_k, y_1, \ldots, y_k\}$).

Now, let's look at the sum again (no longer fixing $H_P, H_R$, and
$H_{PR}$).  We see that for any given choice of pairs $(x_1, y_1),
\ldots, (x_k, y_k)$ with $\ppi_X(x_i) = \ppi_Y(y_i)^*$, the summands
that involve deleting these pairs and gluing their endpoints are the
same in number for each choice of height assignment for the deleted
pairs.  The coefficient for each height is just
$\frac{\h^k}{(|X|+|Y|)!}$ times the number of height assignments
$H_{PR}$ that restrict to the given height assignment, and also
satisfy $H_{PR}(x_i) > H_{PR}(y_i)$ for all $1 \leq i \leq k$.  In
other words, this is $\h^k$ times the probability of picking a height
assignment randomly of $PR$ that has $x_i$ greater in height than
$y_i$ for all $i$, and is identical with the given height assignment
on all $x, y \notin \{x_1, \ldots, x_n, y_1, \ldots, y_n\}$.  So we
get that the coefficient is $\frac{\h^k}{2^k (|X|+|Y|-2k)!}$. 

But this is exactly what we would expect, desiring that \eqref{ltsiso}
hold. That is because the left-hand side, as described previously in
our discussion of $\frac{\pi^N}{2^N}$, just involves summing over all
$N$ of $\frac{\h^k}{2^k}$ times $\Phi_W$ of the collection of
necklaces described for each choice of pairs $(x_1, y_1), \ldots,
(x_N, y_N)$ with some choice of sign; and then $\Phi_W$ just sums over
$\frac{1}{(|X|+|Y|-2N)!}$ times each choice of height assignment for
this collection of necklaces. The sign choice just matches exactly with
the sign $\prod_{i} [\ppi_X(x_i), \ppi_Y(y_i)]$ appearing in \eqref{hpre},
since each commutator is $-1$ just in the case that $\ppi_Y(y_i) \in Q$.

This proves \eqref{ltsiso} and hence that $\Phi_W$ is an isomorphism
of $\k[\h]$-algebras, using $*_\h$ as the ring structure on $\Sym L[\h]_\moy$.

\subsection{Associativity of $*_\h$.} \label{ass}
Although we already know from commutativity of the diagram and associativity
of $A_\h$ that $*_\h$ is associative, it is easy to prove directly. We
prove
\begin{equation} \label{tpass}
(P *_\h R) *_\h S = P *_\h (R *_\h S)
\end{equation}
where $P, R$, and $S$ are collections of necklaces of the form \eqref{pform}
(with different indices).

First we describe the left-hand side of \eqref{tpass}
Let $X, Y$, and $Z$ be the sets of abstract edges
of $P, R,$ and $S$, and let $\ppi_X: X \rightarrow \dq, \ppi_Y: Y
\rightarrow \dq$, and $\ppi_Z: Z \rightarrow \dq$ be the projections
from occurrences of edges to edges of $\dq$.

We sum over all sets of pairs $\{(x_1, y_1), \ldots, (x_N, y_N)\}
 \subset X \times Y$, such that $y_i = x_i^*$ for each $i$ (and we
 assume that the $x_i$ and the $y_i$ are all distinct).  Summing over
 $\frac{\h^N}{2^N}$ times the necklaces we get by cutting out these
 pairs of edges and gluing their endpoints, we get $P *_\h R$ as
 described in the previous section.

To get $(P *_\h R) *_\h S$, we will first be summing over choices of
pairs $\{(x_1, y_1), \ldots (x_N, y_N)\}$, and then over pairs
$\{(w_1, z_1), \ldots, (w_M, z_M)\} \subset W \times Z$, where $W = (X
\setminus \{x_1, x_2, \ldots, x_N\}) \sqcup (Y \setminus \{y_1, y_2,
\ldots, y_N\})$, and performing a similar operation.  We can also
describe this as summing over pairs $(x_1, y_1), \ldots, (x_N, y_N),
(x_{N+1}, z_1), (x_{N+2}, z_2), \ldots, (x_{N+M_1}, z_{M_1}),$ \\
$(y_{N+1}, z_{M_1+1}), (y_{N+2}, z_{M_1+2}), \ldots, (y_{N+M_2},
z_{M_1+M_2})$, again where all $x_i, y_i,$ and $z_i$ are distinct, and
in each pair, one edge is the reverse of the other.  This description,
along with signs and coefficients, is exactly the same as what we
could obtain in the same way from $P *_\h (R *_\h S)$, proving
associativity.

\section{The Moyal coproduct} 
\subsection{Definition of $\Delta_\h$.}\label{cps}
There is a nice formula for the coproduct on $\Sym L[\h]_\moy$ compatible
with the the $*_\h$ product.  The formula is actually surprisingly
similar to the one for $*_\h$.  We will be giving the coproduct which
makes the diagram \eqref{2d} consist of coalgebra homomorphisms
(namely, the maps $\tr$ and $\Phi_W$ involving $\Sym L[\h]_\moy$); the map
$\Phi_W$ will then be an isomorphism of bialgebras.  The coproduct can
be described as follows: We need to define an operator $e^{\frac{1}{2}
\h \pi'}: \Sym L[\h]_\moy \rightarrow \Sym L[\h]_\moy \o \Sym L[\h]_\moy$.  To do
this, we set
\begin{equation}
\pi' = \sum_{e \in Q} \frac{\partial}{\partial e} \frac{\partial}{\partial e^*}
\end{equation}
and we define operators
\begin{multline} \label{bco}
\frac{\partial}{\partial e_1} \frac{\partial}{\partial e_1^*}
 \frac{\partial}{\partial e_2} \frac{\partial}{\partial e_2^*} \cdots
 \frac{\partial}{\partial e_N}\frac{\partial}{\partial e_N^*}: \\ \Sym
 L[\h]_\moy \rightarrow 
\Sym L[\h]_\moy \o \Sym L[\h]_\moy.
\end{multline}
\textbf{Note that order is important: We expand the exponential keeping
track order, so that the adjacent terms $\frac{\partial}{\partial e_1}
\frac{\partial}{\partial e_1^*}$ have a meaningful relationship.}
Specifically, the operator \eqref{bco} acts as follows: Taking a collection of
necklaces $P = P_1 \& P_2 \& \cdots \& P_n$ \\ $\& 
v_1 \& v_2 \& \cdots \& v_q \in \Sym L[\h]_\moy$, where
each $P_i \in L$ is a cyclic monomial (i.e.~a necklace), let $X$ be
the set of abstract edges of $P$ and $\ppi_X: X \rightarrow \dq$ the
projection (cf.  Section \ref{mps}).  Then we sum over all choices of
pairs $(x_1, y_1), \ldots, (x_N, y_N)$ such that the $x_i$ and $y_i$
are all distinct (the set $\{x_1, y_1, \ldots, x_N, y_N\}$ has $2N$
elements), and $\ppi_X(x_i) = \ppi_X(y_i)^*$ for all $i$. We delete the
edges and glue the endpoints, obtaining another collection of
necklaces.  More precisely, the cutting and gluing is done as
described in Section \ref{mpds}, except that the pairs of abstract edges
are pairs of elements of the same set $X$ (there is only one collection
of necklaces). Here, $I = \{x_1, y_1, x_2, y_2,
\ldots, x_N, y_N\}$ and $\phi(x_i) = y_i$ for all $i$.  Now, the only
difficult part is figuring out what components to assign to each
necklace (the first or second), and what sign to attach to each
choice.

We sum over all component assignments of the resulting chain of
necklaces: suppose that the above procedure yields the collection $R_1
\& R_2 \& \cdots \& R_m \in \Sym L[\h]_\moy$ (this includes the original
idempotents $v_1, v_2, \ldots, v_q$); then the contribution to the
result of \eqref{bco} applied to $P$ is the following:

\begin{equation} \label{cass}
\sum_{\mathbf{c} \in \{1,2\}^m} s(\mathbf c, I, \phi)
 R_1^{c_1} \& R_2^{c_2} \& \cdots \& R_m^{c_m},
\end{equation}
where $R_i^{c_i}$ denotes $R_i \o 1$ if $c_i = 1$ and $1 \o R_i$ if
$c_i = 2$, and the symmetric product in $\Sym L[\h]_\moy \o \Sym L[\h]_\moy$
is the expected $(X \o Y) \& (X' \o Y') = (X \& X') \o (Y \& Y')$,
with $1 \& X = X \& 1 = X$ for all $X$.  The term $s(\mathbf c, I, \phi)$ 
is a sign which is determined as follows: 
\begin{equation}
s(\mathbf c, I, \phi) = s_1 s_2 \cdots s_n,
\end{equation}
where $s_i = 1$ if the component assigned to the start of
arrow $x_i$ is $1$ and the component assigned to the target of arrow
$x_i$ is $2$; $s_i = -1$ if the component assigned to the start of
arrow $x_i$ is $2$ and the component assigned to the target of arrow
$x_i$ is $1$; and $s_i = 0$ if the start and target are assigned the
same component.  

Let's more precisely define what it means to say ``the component
assigned to the start/target of an arrow'' which is deleted from $P$.
What we mean by this is simply the orbit of the arrow $x_i$ in $X$
under $f$.  Each orbit corresponds to one of the $R_i$.  So, there is
a map $g: X \rightarrow \{1, 2, \ldots, m\}$, which corresponds to
which $R_i$ the ``start'' of each edge is assigned to.  The
``targets'' are the same as the ``starts'' of the next edge, so that
$g(x+1)$ gives the component that the ``target'' of $x$ is assigned
to. Here the ``$+1$'' operation is once again $(i,j)+1=(i,j+1)$ mod
$l_i$, or in other words, $x+1$ is the edge succeeding $x$.

We then have that
\begin{equation} \label{scc}
s_i = \begin{cases} 1 & c_{g(x_i)} < c_{g(x_i)+1}, \\
                    0 & c_{g(x_i)} = c_{g(x_i)+1}, \\
                   -1 & c_{g(x_i)} > c_{g(x_i)+1}.
\end{cases}
\end{equation}

This assignment of signs has a combinatorial flavor because it is
essentially what the ``colorings'' of \cite{S} reduce to.  There does
not seem to be a way to avoid this complication in choosing signs,
because the sign is what prevents the coproduct from being
cocommutative.

As before, we extend linearly to powers $(\pi')^N$ and to $e^{\frac{1}{2} \h \pi'}$.
Then, the coproduct is given by
\begin{equation}
\Delta_\h = e^{\frac{1}{2} \h \pi'}: \Sym L[\h]_\moy \rightarrow \Sym L[\h]_\moy \o \
\Sym L[\h]_\moy,
\end{equation}
and as before we can describe this action on our element $P$ as
\begin{equation}
\Delta_\h(P) = \sum_{(I, \phi)} \frac{\h^{\#(I)/2}}{2^{\#(I)/2}} P_{I, \phi}, \\ 
\end{equation}
where the sum is over all $I \subset X$ with involution $\phi$ such
that $\ppi_X \circ \phi = * \circ \ppi_X$, and the element $P_{I, \phi}$
is given from the result of the cuttings and gluings by summing over
component assignments as described in \eqref{cass}.

\subsection{Proof that $\Delta_\h$ is obtained from $\Phi_W$.} \label{cpts}
Let's prove that this coproduct $\Delta_\h$ makes the diagram \eqref{2d}
consist of coalgebra homomorphisms.  It suffices to prove that
$\Phi_W$ is a coalgebra homomorphism.

Take an element $P$ of the form \eqref{pform} with set of abstract
edges $X$ and projection $\ppi_X: X \rightarrow \dq$.  Now, let's
consider what the element $\Delta(\Phi_W(P))$ is in $A$. We know that
for each height assignment $H_P$ of $P$, $\Delta(P_{H_P})$ involves
summing over all pairs $(I, \phi)$ with $I \subset X$ and $\phi: I
\rightarrow I$ an involution satisfying $\ppi_X \circ \phi = * \circ
\ppi_X$, cutting and gluing as before. Then we sum over all component
assignments such that if $x, y \in I$ with $\phi(x) = y$, and the
heights satisfy $H(x) < H(y)$, then the component assigned to the
start of $x$ is $1$ and the component assigned to the target of $x$ is
$2$. When the components cannot be assigned in this way, this pair
$(I, \phi)$ cannot be used.  These notions are all explained more
precisely in the preceding section.

Then we multiply in a sign $s(I, \phi, H)$ and a power of $\h$ determined
as follows: for each pair $x, y \in I$ with $\phi(x) = y, H(x) <
H(y)$, we multiply a $+1$ if $x \in Q, y \in Q^*$ and a $-1$ if $x \in
Q^*, y \in Q$.  We also multiply in $\h^{\#(I)/2}$ (note: this power
of $\h$ is different from the one in \cite{S} simply because we are
describing the structure for $A_\h$, not $A$: it is easy to see in
general how the relations for the algebra and the formula for
coproduct change if we introduce a new formal parameter $\hbar$ into 
(3.3) of \cite{S}).

So we find that $\Delta(P_H)$ is just a sum over cuttings and gluings,
and over component choices $\mathbf c$ compatible with the heights;
our sign choice satisfies $s(I, \phi, H) = s(\mathbf c, I, \phi)$,
where $I = \{x_1, y_2, \ldots, x_m, y_m\}$, and for all $i$, $x_i \in
Q$ and $\phi(x_i) = y_i$; finally, we multiply in $\h^m$ for cuttings
and gluings involving $\#(I)=2m$.

Hence, $\Delta(\Phi_W(P))$ is just given by a sum over all cuttings and
gluings $(I, \phi)$ together with component choice $\mathbf c$, 
multiplying in $\h^{\#(I)/2}$, the sign $s(\mathbf c, I, \phi)$, and
the coefficient $\frac{1}{\#(P)!}$ where $\#(P)$ is the number of edges
in $P$, i.e.~the total number of height assignments. 

Each summand in $\Delta(\Phi_W(P))$ is clearly given by a height
assignment of the term in $\Delta_\h(P)$ corresponding to the same $(I,
\phi, \mathbf c)$.  For each term in $\Delta_\h(P)$, the coefficients of
the height-assigned terms in $\Delta(\Phi_W(P))$ are all the same.  So
we see that $\Delta(\Phi_W(P)) = (\Phi_W \o \Phi_W)(P')$, for some $P'
\in \Sym L[\h]_\moy \o \Sym L[\h]_\moy$, where $\Phi_W \o \Phi_W (P \o R) =
\Phi_W(P) \o \Phi_W(R)$.

The element $P'$ can be computed just as we were computing $\Delta(\Phi_W(P))$, but instead of multiplying in $\frac{1}{\#(P)!}$, we need to multiply by the 
fraction of all height choices compatible with this
component choice.  But clearly, each pair $x, y \in I, \phi(x) = y$
induces a single restriction on the choice of heights, namely that
$H(x) < H(y)$ if the component assigned to the start of $x$ is $1$
and the component assigned to the target of $x$ is $2$, and $H(y) > H(x)$
if the opposite is true (the start of $x$ is assigned component $2$ and
the target assigned $1$).  Note that the component assigned to the start
and target of $x$ cannot be the same for there to exist any compatible
height choices.

We see then that, provided a compatible height choice exists, we have
$\#(I)/2$ restrictions, each of which occur with independent probabilities
$\frac{1}{2}$. Hence the coefficient is just $\frac{1}{2^{\#(I)/2}}$. This
shows that $P' = \Delta_\h(P)$, proving that $\Phi_W$ is a coalgebra homomorphism
and hence an isomorphism of bialgebras. (In fact we have now proved that
$(\Sym L[\h]_\moy, *_\h, \Delta_\h)$ is a bialgebra, since we have
transported the multiplication and comultiplication from the bialgebra
from \cite{S}.  But, it is possible to give a direct proof that this is
a bialgebra, which we omit in this version: see \cite{GS}).

\subsection{Coassociativity of $\Delta_\h$.} \label{casss}
Using the coassociativity of $A_\h$ from \cite{S}, we already know from
the fact that $\Phi_W$ is an isomorphism that the product $\Delta_\h$ is
coassociative, but it is not difficult to prove directly. We do that here
by proving
\begin{equation} \label{coasse}
(1 \o \Delta_\h) \Delta_\h (P) = (\Delta_\h \o 1) \Delta_\h (P),
\end{equation}
where once again $P$ is of the form \eqref{pform}.

The left-hand side can be described by summing over choices of cutting
pairs and components $(I, \phi, \mathbf c)$ for $P$, and then cutting
pairs and components for the first component of the result, $(I',
\phi', \mathbf {c'})$, and gluing, assigning the components, etc., and
multiplying by a sign and power of $\frac{\h}{2}$.  We see that this
is the same as choosing just once the triple $(I'', \phi'', \mathbf
{c''})$, where $\mathbf {c''}$ assigns each necklace to one of three
components, $1, 2,$ or $3$, $I'' = I \cup I'$, and $\phi'' |_I = \phi,
\phi''|_{I'} = \phi'$. Then we can cut and glue just one time to get a
tensor in $\Sym L[\h]_\moy^{\o 3}$; the sign and power of $\frac{\h}{2}$
can be determined by using \eqref{scc} where now the two sides of the
inequality have values in $\{1,2,3\}$.

For the same reason, the right-hand side of \eqref{coasse} can be
described in the preceding way, proving \eqref{coasse} and hence
coassociativity.

\section{The antipode} \label{as}
Using $\Phi_W$ and the formula for the antipode in \cite{S}, it
immediately follows that our antipode $S: \Sym L[\h]_\moy$
is given by the formula
\begin{equation} \label{man}
S(P_1 \& P_2 \& \cdots \& P_m) = (-1)^m P_1 \& P_2 \& \cdots \& P_m,
\end{equation}
where each $P_i \in L$ is a necklace (i.e.~a cyclic monomial or vertex
idempotent).  It is immediate that $S^2 = \Id$.  Indeed, $S$ is
diagonalizable with eigenvalues $\pm 1$ and eigenvectors which are
collections of necklaces of the form \eqref{pform}.

Unfortunately, a direct proof that \eqref{man} is the antipode for $\Sym
L[\h]_\moy$ turned out to be too difficult.  The authors are
interested in any good proof of this fact from purely the Moyal point of
view.

\subsection{Acknowledgements}
The authors would like to thank Pavel Etingof for useful
discussions. Also, thanks to Kevin Costello, Mohammed Abouzaid, and Ben
Wieland for discussing Remark \ref{orr}.  The work of both authors was
partially supported by the NSF.

\bibliography{moyal}
\bibliographystyle{amsalpha}

\footnotesize{
{\bf V.G.}: Department of Mathematics, University of Chicago, 
Chicago IL
60637, USA;\\ 
\hphantom{x}\quad\, {\tt ginzburg@math.uchicago.edu}}
\smallskip 

\footnotesize{
{\bf T.S.}: Department of Mathematics, University of Chicago, 
Chicago IL
60637, USA;\\ 
\hphantom{x}\quad\, {\tt trasched@math.uchicago.edu}}

\end{document}